\documentclass[preprint,12pt]{elsarticle}
\usepackage[utf8]{inputenc}
\usepackage{lmodern}
\usepackage{amsmath, amssymb, amsthm, bm}
\usepackage{graphicx}
\usepackage{booktabs}
\usepackage{hyperref}
\usepackage{float}
 
\theoremstyle{definition}

\newtheorem{definition}{Definition}
 
\newcommand{\Dl}{\bm{\Delta}}
\newcommand{\TT}{\mathbb{T}^{2}}
\newcommand{\TF}{\mathbb{T}^{4}}
\newcommand{\dv}{\Delta v}

\usepackage{xcolor}

\journal{Journal of Computational Physics}
\biboptions{sort&compress}
 
\begin{document}
 
\begin{frontmatter}
 
\title{Low Residuals, Wrong Value: Learning the Viscosity Solution of  Hamilton--Jacobi--Isaacs Equations with PINNs}
 
\author[ub1,icrea]{Raul Jimenez}
\author[ub1,ub2]{Pablo Tejerina-P\'erez\corref{cor1}}
\ead{pablo.tejerina@icc.ub.edu}
\author[ub1,ub2]{Leonid Sarieddine}
\author[ub1,ub2]{Pedro Taranc\'on-\'Alvarez}
\author[ub1,ub2]{Marisol Traforetti}
\author[ub1,ub2]{Pau Sol\'e-Vilar\'o}
\author[ucu]{Oleksii Ignatenko}
\author[kyiv]{Nazar Tkach}
\author[harvard]{Pavlos Protopapas}

\cortext[cor1]{Corresponding author.}

\affiliation[ub1]{organization={Institut de Ciències del Cosmos (ICCUB), Universitat de Barcelona}, city={Barcelona}, country={Spain}}
\affiliation[icrea]{organization={Institució Catalana de Recerca i Estudis Avançats (ICREA)}, city={Barcelona}, country={Spain}}
\affiliation[ub2]{organization={Departament de Física Quàntica i Astrofísica, Universitat de Barcelona}, city={Barcelona}, country={Spain}}
\affiliation[ucu]{organization={Faculty of Applied Sciences, Ukrainian Catholic University}, city={Lviv}, country={Ukraine}}
\affiliation[kyiv]{organization={Department of Computer Science, National University of Kyiv-Mohyla Academy}, city={Kyiv}, country={Ukraine}}
\affiliation[harvard]{organization={Institute for Applied Computational Science, Harvard University}, city={Cambridge, MA}, country={USA}}
 
\begin{abstract}
Physics-informed neural networks (PINNs) and other residual-trained solvers judge a candidate solution by the size of its sampled PDE residual. For first-order Hamilton-Jacobi equations this is not sufficient: the same equation and boundary data admit Lipschitz almost-everywhere solutions that differ from the viscosity solution by $O(1)$ while violating the equation only on sets of measure zero. We study whether residual minimization selects such solutions in a fixed-horizon pursuit-evasion game whose Hamilton-Jacobi-Isaacs value function is known in closed form on the circle $S^1$ and on the flat tori $\mathbb{T}^{2}$ and $\mathbb{T}^{4}$. We train four neural parameterizations with ten seeds on each geometry, 120 runs in total. A standard network with periodic Fourier features converges to folded almost-everywhere solutions in all 30 runs: small residuals and the expected capture-layer gradient, yet values that are wrong over most of the domain. A network with periodicity imposed through the loss recovers the viscosity solution in one and two dimensions but becomes unreliable in four. Two parameterizations that receive the geodesic distance do substantially better: a single network fed the wrapped displacement and its distance recovers the viscosity solution in all 30 runs, while a minimum over periodic image charts succeeds in 27 of 30 and introduces a wrong-image failure mode. Finally, for the convex Hamiltonian considered here, we use the semiconcavity of the viscosity solution to build a mesh-free finite-difference test that needs only forward network evaluations; it detects the folded solutions and the harmful min-of-charts failures without access to the exact value function. A small sampled residual thus does not by itself validate a learned Hamilton-Jacobi solution; the parameterization and independent structural diagnostics can be essential for selecting the viscosity solution.
\end{abstract}
 

\end{frontmatter}

\section{Introduction}
\label{sec:intro}
 
The standard argument for trusting a physics-informed neural network (PINN) solver \citep{Raissi2019, sirignano2018dgm, Karniadakis2021} is that, given a partial differential equation (PDE) with a unique solution, if the network's solution yields a PDE residual that is small on a dense set of collocation points, the network's solution should approach the true solution as the residual vanishes. For first-order Hamilton--Jacobi (HJ) equations this reasoning breaks down. The solution of interest is the so-called viscosity solution~\citep{crandall_lions_viscosity, crandall1984some}; for the Hamilton--Jacobi--Isaacs (HJI) equation of a differential game it is the value function~\citep{Isaacs1965, evans1984differential}, the minimax solution of Krasovskii and Subbotin~\citep{krasovskii_subbotin1988, subbotin1995}. The same equation, however, also admits infinitely many \emph{almost-everywhere
solutions}\footnote{Concrete definitions for viscosity and almost-everywhere solutions are given below in Section \ref{sec:problem statement}.}~\citep{crandall_lions_viscosity, cannarsa2004semiconcave}:
Lipschitz functions that match the boundary data and satisfy the equation
outside a set of measure zero, on which their gradient jumps, but that differ from the value function, often by $O(1)$. These classical zigzag
constructions can make a sampled $L^{2}$ residual small. Whether a residual-trained network actually finds them is an empirical question that
matters for neural HJ solvers in control and reachability~\citep{han2018deepbsde, bansaltomlin2021deepreach,
onken2021neuralhamiltonian}, where HJI equations describe the reachable sets of dynamic games~\citep{mitchell2005time, bansal2017hamilton} and whose interesting regimes offer no ground truth against which a wrong answer could be noticed.

To study whether this ambiguity matters in practice, we consider a fixed-horizon pursuit-evasion game for which the viscosity solution is known in closed form. We study the same game on three translation-invariant geometries: the circle $S^1$, the two-dimensional torus $\TT$, and the four-dimensional torus $\TF$. This allows us to compare the learned function directly with the true value function, while following the same problem from one to four spatial dimensions. For each geometry, we train four parameterizations of the value with the same optimizer, batch sizes, and number of iterations with ten seeds each.

We compare four different neural parameterizations of the value function. A standard network based on periodic Fourier features consistently reaches a small residual but converges to the wrong almost-everywhere solution on every geometry and every seed. By contrast, parameterizations that incorporate more of the geometry perform much better. Enforcing periodicity through the loss is sufficient in one and two dimensions but becomes unreliable in four, while the two parameterizations that receive the geodesic distance recover the viscosity solution accurately in almost all runs.

These experiments show that a small sampled residual is not, by itself, enough to determine which solution of a first-order HJ equation has been learned. They also show that the choice of parameterization can strongly affect solution selection. In the present problem, giving the network the geodesic-distance structure is enough to recover the correct value function.

This leaves a practical question: how can one detect a wrong solution when the exact value function is unavailable? For the convex Hamiltonian considered here, the viscosity solution satisfies a one-sided semiconcavity bound. We use this property to construct a finite-difference test that requires only forward evaluations of the trained network. In our experiments, the test clearly distinguishes the folded solutions from the correct ones, and identifies a new type of failure mode particular to one of the parameterizations whose residual alone failed to highlight.

The paper is organized as follows. Section~\ref{sec:related} reviews previous work. Section~\ref{sec:problem statement} introduces the differential game, its HJI equation and the closed-form value function (the viscosity solution), and shows that the equation has infinitely many almost-everywhere solutions besides it. Section~\ref{sec: methodology} presents the PINN loss, explains why it cannot distinguish the value function from a folded solution, and describes the four network configurations. Section~\ref{sec:results} compares the trained PINNs with the closed-form solution and analyzes the different ways in which the parametrizations succeed or fail, and Section~\ref{sec:conclusions} offers some conclusions and discussions. The training details are given in \ref{app:training}, while the semiconcavity test is presented in \ref{app:cert}.

\section{Previous work}
\label{sec:related}

When it comes to closely related applications of PINNs in game theory, in a parallel work~\citep{companion2026embeddings} we study the same pursuit-evasion game on the flat torus $\TT$ from the side of the players. Both players were represented by multi-head PINNs, trained against each other by self-play, and the learned strategies were compared with the optimal ones, which are known in closed form for this geometry (for other applications of multi-head PINNs, see~\citep{tarancon2025efficient, JimenezEmbeddings2026, taranconalvarez2026neural}). There, the play was also checked at the level of the value function by a separate network trained on the residual of the HJI equation of the game. The residual level of that network turned out not to reflect the quality of its values. This paper studies why, by asking which solution of the HJI equation a residual-trained network converges to.

The question regarding which solution of a first-order HJI equation is the correct one has long been settled. The value function of a differential game is characterized, under suitable assumptions, as the unique viscosity solution of the associated HJI equation~\citep{crandall_lions_viscosity}. Equivalently, it can be characterized as the unique minimax solution in the sense of Krasovskii and Subbotin~\citep{krasovskii_subbotin1988, subbotin1995}. When the Hamiltonian is convex, as in our problem, the viscosity solution is also characterized by a one-sided condition on its creases, namely semiconcavity~\citep{cannarsa2004semiconcave}. This is the HJ counterpart of Oleinik's entropy condition for conservation laws~\citep{oleinik1957}: in one spatial dimension, the derivative of a solution of a HJ equation solves a conservation law, and for a convex Hamiltonian the semiconcavity of the solution corresponds to Oleinik's one-sided condition on its derivative~\citep{cannarsa2004semiconcave}. For HJI equations whose Hamiltonian is neither convex nor concave, no such one-sided condition exists~\citep{cannarsa2004semiconcave}. Characteristics describe the value where it is smooth, and Melikyan's method of singular characteristics extends them across the surfaces on which it is not~\citep{melikyan1998}, including for a differential game of simple approach on manifolds~\citep{melikyan_ovakimyan1993}. Which nonsmooth continuation is admissible is decided by the generalized-solution framework above. Section~\ref{sec:problem statement} defines these notions for our game. This theory says which solution is correct, but not whether a residual-trained network finds it.

Most studies of failure in physics-informed training concern optimization: loss landscapes and curricula~\citep{krishnapriyan2021failure}, gradient pathologies~\citep{wang2021understanding, wang2022when}, trivial solutions that satisfy the residual at sparse collocation points~\citep{leiteritz2021trivial}, and propagation failures addressed by residual-driven resampling~\citep{daw2023r3}. The failure studied here concerns what the training converges to rather than how. The estimates of \citet{mishra2023estimates}, which bound the error of a PINN by its residual, rely on a stability property that first-order HJ equations have only among viscosity solutions. \citet{wang2022l2pinn} reach a related conclusion for Hamilton--Jacobi--Bellman equations of stochastic control, for which a small $L^{2}$ residual need not imply a small error, and train with an $L^{\infty}$ loss instead.

Several works make neural solvers select the viscosity solution. \citet{Darbon2016HJB} compute viscosity solutions of certain convex-Hamiltonian equations from Hopf-type representation formulas, and \citet{darbon2021architectures} build min-plus networks that represent such solutions exactly, without training(see also~\cite{darbon2020overcoming, darbon2023minplus}), in the tradition of idempotent methods~\citep{mceneaney2007curse}. Our min-of-charts network (Section~\ref{sec:charts}) applies the same idea to a trained network on a problem with a capture set, a finite horizon and a periodic geometry, which those representations do not cover. \citet{esteve2024fdls} replace the pointwise residual by that of a monotone finite-difference scheme, which selects the viscosity solution~\citep{crandall1984two, barles1991convergence}, and \citet{bokanowski2026monotone} add error bounds and a multi-level training scheme. \citet{kim2026monotone} report that a continuous residual converges to a wrong solution on a one-dimensional problem, and that a finite-difference residual, combined with refinement near the kink, largely removes it. \citet{golpashin2026viscosity} find that a PINN for a one-dimensional eikonal problem ends at a wrong solution on about half of ten seeds, and enforce the viscosity conditions through a min--max loss over test functions. In general, the above works employ methodologies that change the loss, while in this work the plain residual is kept and instead the parameterization is changed.

Other works check or constrain a solution after or during training. \citet{albert2002posteriori} give a posteriori error estimates for HJ equations that apply to any numerical approximation, and \citet{lin2023assurances} obtain probabilistic guarantees for neural reachability after training. In neural signed-distance fields, second-order regularizers steer training away from spurious minimizers~\citep{benshabat2022digs, wang2023singularhessian}. \citet{deryck2022wpinns} select entropy solutions of conservation laws through an adversarial weak formulation, and placing a known singularity in the parameterization is standard for the point source of eikonal networks~\citep{binwaheed2021pinneik}.

\section{Problem statement}
\label{sec:problem statement}

We consider a pursuit-evasion game on flat, compact geometries. We first introduce the game, the Hamilton--Jacobi--Isaacs (HJI) equation for the value function of the game, and its closed-form viscosity solution. We then explain in what sense a non-differentiable function solves this equation, show that the equation has infinitely many such solutions besides the viscosity solution, and state the questions addressed in this paper.

\subsection{The pursuer and evader game}
\label{subsec:game}

Let us consider a pursuer and an evader moving on the flat torus $\mathbb{T}^{n}$, represented by the cube $[0,L)^{n}$ with opposite faces identified. We take $L = 2$ and study the case $n = 1$, which corresponds to the circle $S^{1}$, and the cases with $n=2,\, 4$ corresponding to the tori $\mathbb{T}^{2}$ and $\mathbb{T}^{4}$, respectively. The pursuer and the evader move at constant speeds $v_P = 1.2$ and $v_E = 1$, respectively, and can freely choose their directions of motion. The evader is caught when their distance drops to $\varepsilon = 0.08$. The game has a finite horizon, $\tau$, which denotes the \emph{remaining horizon}, the time ``still available for play''. This horizon decreases at unit rate as play proceeds, and the game ends either at capture or when $\tau$ reaches $0$, whichever comes first. Here, we consider horizons $0 < \tau \le \tau_{\max} = 8$. By translation invariance, the game depends only on the relative position $\bm{x}_E - \bm{x}_P$, defined modulo $L$ in each coordinate via the wrapping function $\bm{w}(\bm{x})$. We define the state as the relative displacement vector, $\Dl$, in the \emph{fundamental domain} $[-L/2,L/2)^{n} = [-1,1)^{n}$
\begin{equation}
\label{eq:wrap}
\Dl = \bm{w}(\bm{x}_E - \bm{x}_P), \qquad
w_i(\bm{x}) = \bigl[(x_i + L/2) \bmod L\bigr] - L/2,
\end{equation}
 which places the pursuer at the origin. The geodesic distance between the players, or minimum-image distance, is defined as:
\begin{equation}
\label{eq:geodesic}
d = \min_{\bm{s} \in L\mathbb{Z}^{n}} \|\bm{x}_E - \bm{x}_P - \bm{s}\|
  = \bigl\|\bm{w}(\bm{x}_E - \bm{x}_P)\bigr\| = \|\Dl\| ,
\end{equation}
with maximum $\sqrt{n}$ at the corners of the fundamental domain. This is the distance from the evader to the nearest \emph{periodic image} of the pursuer, where these are defined as $\bm{x}_P + \bm{s}$, with $\bm{s} \in L\mathbb{Z}^{n}$. The states where the nearest image is not unique, so that the evader can be reached along several shortest paths, form the \emph{cut locus} of the pursuer, i.e. the faces $\Delta_i = \pm 1$. There, $d$ is continuous but not differentiable.

The \emph{value function} $V(\Dl,\tau)$ is the capture time under optimal play, capped at time $\tau$. This function is minimized by the pursuer, and maximized by the evader. It satisfies, in a sense made precise below, the HJI equation~\citep{Isaacs1965, evans1984differential} 
\begin{equation}
\label{eq:hji}
\partial_\tau V = 1 - \dv\,\|\nabla V\|, \qquad \dv := v_P - v_E = 0.2,
\end{equation}
for $d(\Dl) > \varepsilon$ and $0 < \tau \le \tau_{\max}$, where $\nabla$ denotes the gradient with respect to $\Dl$. The system is constrained by the terminal and capture conditions: 
\begin{align}
\label{eq:terminal}
&V(\Dl, 0) = 0 \quad\text{ for } \Dl\in\mathbb{T}^n,\\
\label{eq:capture}
&V(\Dl, \tau) = 0 \quad \text{ for } \Dl\in\mathcal{C} \text{ and } 0<\tau\leq\tau_{\rm max},
\end{align}
where $\mathcal{C}=\{\Dl\in\mathbb{T}^n: d(\Dl)\leq\varepsilon\}$ is the \emph{capture set}.
The gradient term in Eq. \eqref{eq:hji} combines the evader's contribution $+v_E\|\nabla V\|$ and the pursuer's one $-v_P\|\nabla V\|$. In Hamiltonian form, the  HJI equation reads
\begin{equation}
\label{eq: HJI in hamiltonian form}
\partial_\tau V + H(\nabla V) = 0, \qquad H(\bm{p}) = \dv\,\|\bm{p}\| - 1 ,
\end{equation}
where the Hamiltonian $H$ is convex since $v_P > v_E$. This holds only in the relative position $\Dl$. In the players' own coordinates, the evader contributes a concave term and the Hamiltonian is neither convex nor concave.

In this setup, the value function, which we identify below as the viscosity solution of \eqref{eq: HJI in hamiltonian form}, can be written as follows:
\begin{equation}
\label{eq:vstar}
V^{*}(\Dl,\tau) = \min\left(\tau,\
\frac{\left(d(\Dl) - \varepsilon\right)_{+}}{\dv}\right),
\qquad (x)_{+} := \max(x, 0).
\end{equation}
From this, we can distinguish two different situations. On the \emph{distance-limited} branch, $(d - \varepsilon)/\dv \le \tau$, the pursuer runs towards the evader along a shortest path, the evader runs directly away, and the gap $d - \varepsilon$ closes at rate $\dv$, so $V^{*} = (d - \varepsilon)/\dv$ and $\|\nabla V^{*}\| = 1/\dv$. On the \emph{horizon-limited} branch the evader cannot be caught in time, so $V^{*} = \tau$ and $\nabla V^{*} = 0$. In every dimension, $V^{*}$ depends on $\Dl$ only through $d$.

\begin{figure}[!tbp]
\centering
\includegraphics[width=\textwidth]{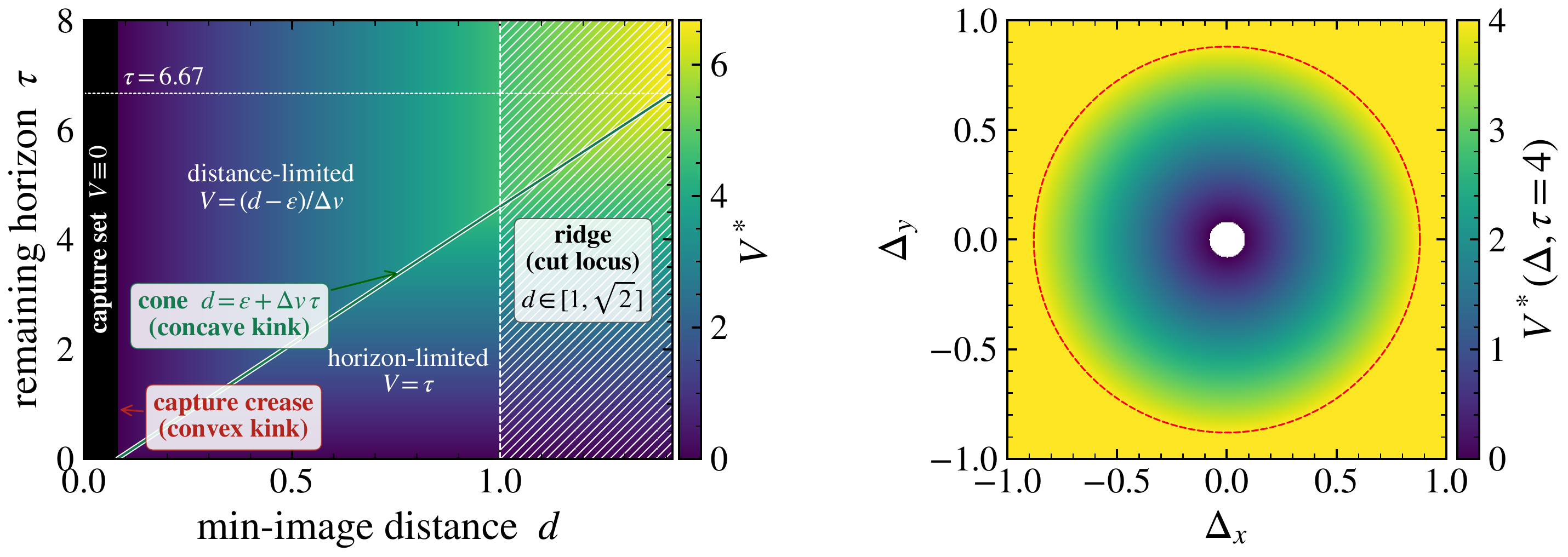}
\caption{The value function $V^{*}$ of Eq. \eqref{eq:vstar} for $\varepsilon = 0.08$, $\dv = 0.2$. \textbf{Left:} $V^{*}$ on $\TT$ in the $(d,\tau)$ plane. We show the capture set (black band), the cone $d_c = \varepsilon + \dv\,\tau$ (green line), the capture crease $d = \varepsilon$ (a convex kink on the boundary), the horizon at which the cone reaches the point of $\TT$ farthest from the pursuer (dotted line at $\tau=6.67$), the nearest points of the cut locus (dashed $d=1$) and the distances spanned by the cut locus (hatched band). \textbf{Right: } A time-slice $V^*(\Dl,\tau=4)$ on the fundamental domain of $\TT$. The dashed, red circle marks the cone $d_c=\varepsilon+\Delta v\,\tau$; inside it the value function grows linearly with $d$, while outside it is $V^*=\tau$.}
\label{fig:vstar}
\end{figure}

Figure~\ref{fig:vstar} (Left) shows $V^{*}$ on $\TT$ in the $(d,\tau)$ plane. The line $d_c = \varepsilon + \dv\,\tau$, where the time to capture equals the remaining horizon, separates the distance-limited branch (upper left) from the horizon-limited one (lower right). We call it the \emph{cone}, since in $(\Dl,\tau)$ space it is a cone around the pursuer, truncated by the faces of the fundamental domain once $\varepsilon + \dv\,\tau > 1$. The hatched band shows the distances $1 \le d \le \sqrt{2}$ spanned by the cut locus.

Figure~\ref{fig:vstar} (Right) shows $V^{*}$ on $\TT$ as a function of the state $\Dl = (\Delta_x, \Delta_y)$ at the fixed horizon $\tau = 4$. The dashed red circle is the cone $d_c = \varepsilon + \dv\,\tau$. Inside it, the state is distance-limited, and $V^{*}$ grows linearly with the distance, with slope $1/\dv$, from zero at the edge of the capture set to $\tau$ at the cone. Outside it, the state is horizon-limited, and $V^{*}$ is constant and equal to $\tau$.

The viscosity solution $V^{*}$ is continuous, but its gradient jumps across some surfaces, called \emph{kinks}. A kink is \emph{concave} if the function forms a peak across it, like $-|x|$ at $x = 0$, and \emph{convex} if it forms a valley, like $|x|$. Inside the domain $d(\Dl) > \varepsilon$, $0 < \tau \le \tau_{\max}$ of Eq. \eqref{eq:hji}, $V^{*}$ has two kinks that are both concave. The first one is the \emph{cone} (dashed red circle in Fig. \ref{fig:vstar} (Right)), while the second one is the \emph{ridge}, which is the peak that $d$ has on the cut locus, and that $V^{*}$ inherits where the cut locus is distance-limited. As the cut locus lies on the faces of the fundamental domain, the ridge is invisible in plots on $[-1,1)^{n}$. Figure~\ref{fig:selector} (Right) shows it for $d$, with the pursuer drawn at the corners.

Because of its kinks, $V^{*}$ is not a classical solution of Eq. \eqref{eq:hji}. Nevertheless, we know that a function is Lipschitz continuous if there exists a finite, constant $M>0$ such that $|V(\bm{y}) - V(\bm{y}')| \le M\,\|\bm{y} - \bm{y}'\|$, where $\bm{y} = (\Dl,\tau)$. By Rademacher's theorem, every Lipschitz function is differentiable almost everywhere, which suggests a weaker notion of solution.
\begin{definition}[Almost-everywhere solution]
\label{def:ae}
A function $V:\mathbb{T}^n\times[0,\tau_{\max}]\to\mathbb{R}$ is called an \emph{almost-everywhere solution} of the system (\ref{eq:hji}--\ref{eq:capture}) if it satisfies the following conditions:
\begin{enumerate}
    \item $V$ is Lipschitz continuous.
    \item It satisfies the HJI equation \eqref{eq:hji} at almost every point with $\{d(\Dl) > \varepsilon,0 < \tau \leq \tau_{\max}\}$.
    \item The terminal \eqref{eq:terminal} and capture \eqref{eq:capture} conditions hold.
\end{enumerate}

\end{definition}
\noindent The viscosity solution $V^{*}$ is one almost-everywhere solution, but not the only one. Suppose the evader on $S^{1}$ were also caught at a point $\Dl_0$ lying outside the capture set. Then, the value of the game would be
\begin{equation}
\label{eq:VZ}
V_{Z} = \min\Bigl(V^{*},\ \frac{\operatorname{dist}(\Dl, Z)}{\dv}\Bigr),
\qquad
\operatorname{dist}(\Dl, Z) = \min_{\Dl' \in Z} \bigl\|\bm{w}(\Dl - \Dl')\bigr\| ,
\end{equation}
with $Z = \{\Dl_0\}$. Like $V^{*}$, $V_{Z}$ is a minimum of functions that satisfy Eq. \eqref{eq:hji} almost everywhere, and it satisfies the conditions \eqref{eq:terminal} and \eqref{eq:capture}, so $V_{Z}$ is an almost-everywhere solution. 
However, it is not the value function: at $\Dl_0$ it drops to zero, forming a valley (a convex kink) where $V^{*}$ has none. Any closed set $Z$ of measure zero outside the capture set, such as a curve on $\TT$, gives another one, so there are infinitely many almost-everywhere solutions that are not the value function. We call such functions \emph{folded solutions}. These present kinks that are not present in the value function. We call such kinks \emph{folds}.

\begin{definition}[Viscosity solution~\citep{crandall_lions_viscosity, crandall1984some}]
\label{def:viscosity}
A continuous function $V$ satisfying the conditions \eqref{eq:terminal} and \eqref{eq:capture}  is a \emph{viscosity solution} of the HJI equation \eqref{eq: HJI in hamiltonian form} if it is both a sub- and a supersolution in the following sense.
For every $C^{1}$ function $\varphi$ and every point $(\Dl,\tau)$ with $d(\Dl) > \varepsilon$ and $0 < \tau < \tau_{\max}$: if $V - \varphi$ has a local maximum, then $\partial_\tau \varphi + H(\nabla \varphi) \le 0$ ($V$ is a \emph{subsolution}); if $V - \varphi$ has a local minimum, then $\partial_\tau \varphi + H(\nabla \varphi) \ge 0$ ($V$ is a \emph{supersolution}).

\end{definition}
\noindent At points where $V$ is differentiable, both conditions reduce to \eqref{eq:hji}. At a peak, $\varphi$ can only touch $V$ from above, with derivatives that are convex combinations of those on either side, and the convexity of $H$ makes the subsolution condition hold. 
At a valley, $\varphi$ can only touch $V$ from below, and the supersolution condition can fail. This occurs for $V_{Z}$: the constant test function $\varphi \equiv 0$ touches it from below at $\Dl_0$. Since $\partial_\tau\varphi=0$ and $\nabla_{\Dl}\varphi=0$, the supersolution inequality gives $H(0) = -1 < 0$.

Semiconcavity makes precise the idea that all kinks are peaks. A function $V$ is \emph{semiconcave} if, on every compact subset of the domain, there is a constant $C$ such that
\begin{equation}
\label{eq:semiconcave}
V(\bm{y} + \bm{h}) - 2V(\bm{y}) + V(\bm{y} - \bm{h}) \le C\,\|\bm{h}\|^{2}, \qquad \bm{y} = (\Dl,\tau),
\end{equation}
whenever the segment from $\bm{y} - \bm{h}$ to $\bm{y} + \bm{h}$ lies in that subset. This bounds the curvature from above only. A peak satisfies \eqref{eq:semiconcave}, whereas a valley does not: for $V(x)=|x|$ at $x = 0$, the left-hand side of \eqref{eq:semiconcave} is $2|h| > C h^{2}$ for small $h$. 
For a convex Hamiltonian, semiconcavity singles out the viscosity solution among the almost-everywhere solutions~\citep{cannarsa2004semiconcave}. Since the viscosity solution is unique~\citep{crandall_lions_viscosity, bardi1997optimal} and coincides with the value function~\citep{evans1984differential, subbotin1995, bardi1997optimal}, semiconcavity helps to distinguish the correct solution from the other almost-everywhere solutions. $V^{*}$ is semiconcave, because it is the minimum of the finitely many smooth functions $\tau$ and $(\|\Dl - \bm{s}\| - \varepsilon)/\dv$, with $\bm{s} \in L\{-1,0,1\}^{n}$, so it is the viscosity solution, in agreement with Eq. \eqref{eq:vstar}. Every other almost-everywhere solution fails to be semiconcave somewhere, typically at a valley.

We address the problem of computing $V^{*}$ on $S^{1}$, $\TT$ and $\TF$ with a physics-informed neural network (PINN), which is trained by minimizing the residual of the HJI equation \eqref{eq:hji} at randomly sampled points (see Section~\ref{sec: methodology}). Every almost-everywhere solution has zero residual outside a set of measure zero, so a small residual does not distinguish the viscosity solution $V^{*}$ from a folded solution. We ask:
\begin{enumerate}
\item[(Q1)] Does a residual-trained network converge to $V^{*}$, or to a different almost-everywhere solution such as a folded one?
\item[(Q2)] How does the answer depend on the way the geometry of the torus (for instance its periodicity and its distance) enters the network, and on the dimension~$n$?
\end{enumerate}

\section{Methodology}
\label{sec: methodology}

We approximate the value function with a PINN~\citep{lagaris1998, Raissi2019, Karniadakis2021}. A PINN represents the unknown function by a neural network and trains it to satisfy the differential equation. It does so by minimising the residual of the equation at randomly sampled points~\citep{sirignano2018dgm}, together with the mismatch in the boundary conditions. The derivatives in the residual are computed by automatic differentiation~\citep{baydin2018autodiff}, and no solution data are needed. We first write the loss for our problem, then explain why minimising it need not select the value function, and finally present the four network architectures that we compare.

For a candidate $V$, the residual of Eq. \eqref{eq:hji} is
\begin{equation}
\label{eq: residual HJI}
r[V] = \partial_\tau V - 1 + \dv\,\|\nabla V\| ,
\end{equation}
which vanishes wherever $V$ satisfies the HJI equation \eqref{eq:hji}. The PINN minimises
\begin{equation}
\label{eq:loss}
\begin{split}
\widehat{L}[V] ={}& \frac{1}{N_{\mathrm{int}}}\sum_{i} r[V](\Dl_i,\tau_i)^{2} \;+\; \lambda\,\frac{1}{N_{0}}\sum_{j} V(\Dl_j,0)^{2} \\
&\;+\; \lambda\,\frac{1}{N_{\varepsilon}}\sum_{k} V(\Dl_k,\tau_k)^{2},
\end{split}
\end{equation}
where the $N_{\mathrm{int}}$ interior points $(\Dl_i,\tau_i)$ are drawn uniformly outside the capture set, the $N_{0}$ points $(\Dl_j,0)$ on the terminal slice $\tau = 0$, and the $N_{\varepsilon}$ points $(\Dl_k,\tau_k)$ inside the capture set $\mathcal{C}$. The last two terms in \eqref{eq:loss} impose the boundary conditions \eqref{eq:terminal} and \eqref{eq:capture}, with a weight $\lambda$ relative to the residual. All points are redrawn at every training step. \ref{app:training} gives the sampling ranges, batch sizes, value of $\lambda$, and other details of the training procedure.

\subsection{Loss blindness to almost-everywhere solutions}
\label{sec:analysis}

The loss \eqref{eq:loss} cannot tell the value function from a folded solution. Both satisfy the HJI equation \eqref{eq:hji} everywhere except on their kinks, which are a set of measure zero that random points hit with probability zero, so both have zero loss. A network is smooth and can only approximate a kink by rounding it over some width $\ell$. 
At a fold or a ridge where two sides meet with gradients pointing in different directions, a smooth network rounds the gradient discontinuity across a transition layer of width $\ell$, while the gradients on either sides are generally smaller than $1/\dv$, causing an order one violation of the HJI equation inside this band. The probability of a sampled point falling within this layer is $O(\ell)$, consequently, each fold adds to the loss an amount $O(\ell)$. As the network sharpens the folds, this contribution vanishes, even though the values on either side can remain incorrect by $O(1)$. Training can therefore lower the loss by sharpening the folds of a wrong solution, rather than removing them.

The cone is different. Away from the cut locus, it can be smoothed without any residual cost. On any region where $d(\Dl)$ is differentiable, the soft minimum is
\begin{equation*}
V_{\mu} = -\mu\log\bigl(e^{-\tau/\mu} + e^{-(d-\varepsilon)/(\dv\,\mu)}\bigr),
\end{equation*}
which tends to $V^{*}$ as $\mu \to 0$, and it satisfies \eqref{eq:hji} exactly as long as $\mu>0$. Indeed, $\partial_\tau V_{\mu} = p$ and $\dv\,\|\nabla V_{\mu}\| = 1 - p$, where $p \in (0,1)$ is the weight of the branch $\tau$ in this soft minimum. A network can therefore smooth the cone without any residual cost, unlike a fold or the ridge.

When the loss fails, the error has a definite sign: a folded solution never lies above $V^{*}$. Intuitively, a valley of a folded solution acts like an extra capture set. For instance, $V_{Z}$ in \eqref{eq:VZ} drops to zero at $\Dl_0$ as if the evader could also be caught there, and an extra way of being caught can only shorten the evader's survival. Formally, since $H$ is convex, every almost-everywhere solution is a viscosity subsolution~\citep[Remark~5.3.2]{cannarsa2004semiconcave}, and a subsolution never exceeds the viscosity solution with the same boundary conditions~\citep{crandall_lions_viscosity, bardi1997optimal}. Whether training actually finds folded solutions is an empirical question, which we address in Section~\ref{sec:results}.

\subsection{Network configurations}
\label{subsec: PINN configurations}

All four configurations use the same network $N$, a fully connected $\tanh$ network that receives the time through fixed features of $t = \tau/\tau_{\max}$ (\ref{app:training}). They differ in how the state $\Dl$ enters. The Fourier-feature and soft-periodic NNs only know that the domain is periodic, while the min-of-charts and geometric NNs also receive the geodesic distance. 

\subsubsection{The Fourier-feature NN}
\label{sec:smooth}

The simplest way to make a network periodic is to feed it periodic functions of the state~\citep{dong2021periodic}. The Fourier-feature NN~\citep{tancik2020fourier} uses $K = 8$ harmonics per coordinate,
\begin{equation}
\phi_{\mathrm{per}}(\Dl) = \bigl(\sin k\omega\Delta_i,\ \cos k\omega\Delta_i\bigr)_{i=1,\dots,n,\ k=1,\dots,K} \ \ \ \ ; \qquad \omega = \dfrac{2\pi}{L} \, ,
\end{equation}
and
\begin{equation}
    V(\Dl,\tau) = N(\phi_{\mathrm{per}}(\Dl), t)\,.
    \label{eq: fourier feature NN}
\end{equation}
It is periodic and smooth by construction, so it can only approximate the kinks of $V^{*}$ by smooth functions of these features, and neural networks are known to learn the high-frequency content that such kinks require only slowly~\citep{rahaman2019spectral, wang2021eigenvector}.

\subsubsection{The soft-periodic NN}
\label{sec:rawsoft}

Instead of building periodicity into the network, the soft-periodic NN imposes it through the loss. The network receives the wrapped displacement of Eq.~\eqref{eq:wrap} itself and nothing else,
\begin{equation}
\label{eq:rawsoft}
V(\Dl,\tau) = N\bigl(\Dl,\ t\bigr),
\end{equation}

so its input jumps from $+L/2$ to $-L/2$ across each face of the fundamental domain. The loss \eqref{eq:loss} gets an extra term that penalises any mismatch between the values on opposite faces\footnote{The constraint then holds only up to the value reached by the penalty, as opposed to an exact imposition through the architecture~\citep{lagaris1998, dong2021periodic, lu2021hard, sukumar2022exact}.}
\begin{equation}
\label{eq:seam}
L_{\mathrm{periodic}} = \lambda\,\frac{1}{N_{p}}\sum_{l}\bigl[N(\bm{x}_l^{+},t_l) - N(\bm{x}_l^{-},t_l)\bigr]^{2},
\end{equation}
over $N_p$ random pairs of points $\bm{x}_l^{\pm}$ on opposite faces, with $t_l=\tau_l/\tau_{\rm max}$ (\ref{app:training}). The penalty matches values only, so the gradient may jump across the faces, which is where the ridge of $V^{*}$ lies. Nothing else about the geometry is supplied, and the network has to discover that $V^{*}$ depends on $\Dl$ only through $d$.

\subsubsection{The min-of-charts NN}
\label{sec:charts}

For the min-of-charts NN, we use the fact that the value function has a structure that a network can imitate. The distance-limited part of the analytic solution $V^*$ depends only on the geodesic distance, which is computed as the minimum over the displacements from the different periodic-images of the pursuer to the evader:

\begin{equation}
    \label{eq: Vstar as min of min}
    V^{*}(\Dl,\tau) = \min\left( \tau, \,\min_{\bm{s} \in L\mathbb{Z}^{n}} \frac{(\|\Dl-\bm{s}\| - \varepsilon)_{+}}{\dv}\right).
\end{equation}
The min-of-charts NN builds this minimum over $\bm{s}\in L\mathbb{Z}^n$ into the architecture, as described next.

In the geometries considered in this work, there are only $2^n$ periodic images of the pursuer (duplicates in the real plane corresponding to different wrapping numbers) that can have the nearest image-distance to the evader. On the circle $S^1$, these correspond to the two directions in which the pursuer can go around. In the 2D torus $\TT$, they correspond to four possibilities: the path directly inside the fundamental domain, the path that wraps around the $x$-coordinate, the one wrapping around $y$, or both. Generally in $\mathbb{T}^n$, these possible displacements to the evader can be written as

\begin{equation}
\label{eq: s_m image displacements}
\bm{z}_m = \Dl - \bm{s}_m, \qquad (\bm{s}_m)_i = L\, b_{m,i}\operatorname{sgn}\Delta_i, \qquad \bm{b}_m \in \{0,1\}^{n}.
\end{equation}
Each value of $m=0,\dots, 2^n-1$ defines a region in the torus $\mathbb{T}^n$, called a \emph{chart}.

\noindent Using a network $N$, we parametrize the value function as
\begin{equation}
\label{eq:min of charts NN}
V(\Dl,\tau) = \min_{m} N\bigl(\phi(\bm{z}_m),\ t\bigr),
\end{equation}
where $\phi(\bm{z}) = \bigl(\bm{z},\ \bm{z}^{2},\ \|\bm{z}\|\bigr)$ are the inputs $\phi$--features of the different charts.
The features must be non-periodic, otherwise all charts would compute the same function. One training point costs $2^{n}$ network evaluations. The chart with $\bm{s}_m = \bm{0}$ has displacement $\Dl$ itself, of norm $d$, so the network receives the geodesic distance but is not told which chart carries it. In this way, the network only has to learn the function of a single image, the \emph{chart} chosen by the minimum in \eqref{eq:min of charts NN}.

\begin{figure}[t!]
\centering
\includegraphics[width=\textwidth]{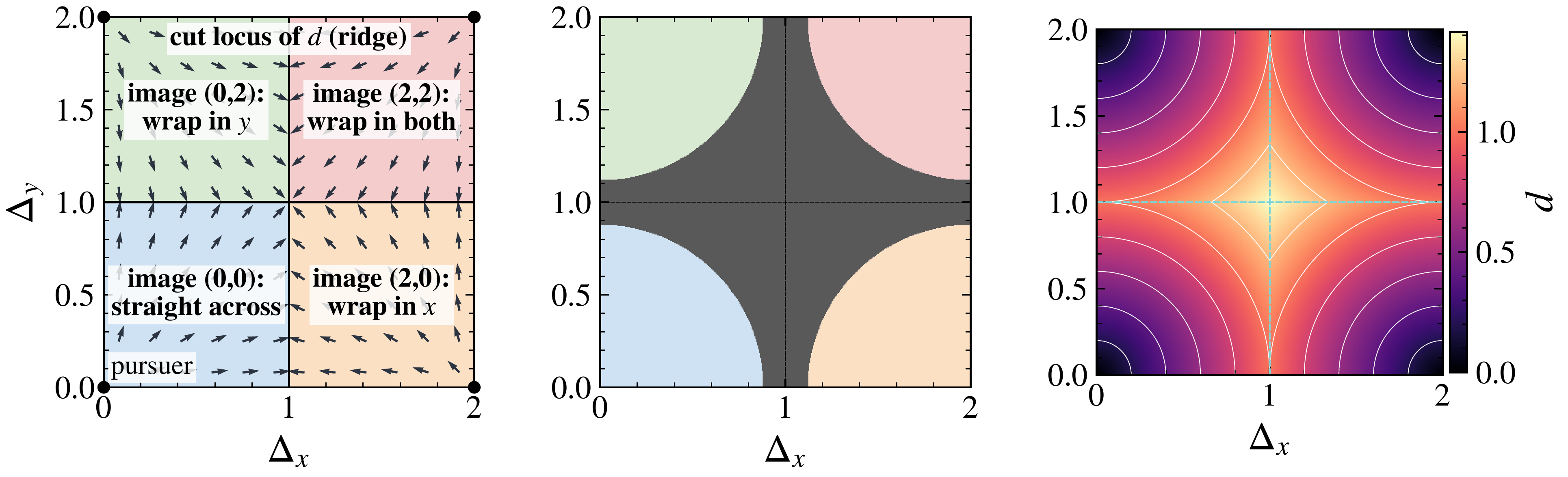}
\caption{The chart selector $m(\Dl,\tau)$ on $\TT$, drawn on $[0,2)^{2}$ with the pursuer at the origin $(0,0)$ and its periodic-images at the other corners. The faces defining the cut locus of the pursuer on the torus, which are $\tilde\Delta_i=\pm 1$ on the fundamental domain, appear as the midlines $\Delta_i = 1$. \textbf{Left.} Analytic nearest-image index $m^{*}(\Dl,\tau)$. The four colors index each of the charts. Arrows show the optimal direction of motion of the pursuer, from the nearest periodic image, towards the evader. The charts meet on the ridge (cut locus). In the fundamental domain $[-1,1)^2$ (coordinates of the text), the nearest image is always in the chart with $\bm{b}_m = 0$. In the representative coordinates of the figure $\Dl\in[0,2)^2$ all signs are positive, so $\bm{z}_m = \Dl - L \bm{b}_m$, the images sit at $L \bm{b}_m$, and the four colors are $\bm{b}_m = \{(0,0), (1,0), (0,1), (1,1)\}$. \textbf{Middle.} Selector $\hat m (\Dl,\tau=4)$ learned by the min-of-charts NN (seed~0). Horizon-limited states, where all charts tie, are shown in grey. On the distance-limited states, the learned $\hat m(\Dl)$ agrees with the analytic value (colors in left plot) everywhere. \textbf{Right.} The nearest-image distance $d =\text{min}_m\,\|\bm{z}_m\|$, equal to $\|\Dl\|$ in the fundamental domain. It is continuous, with a concave kink at the ridge, where the index $m^{*}$ jumps, and $\|\nabla d\| = 1$ on both sides.}
\label{fig:selector}
\end{figure}

This minimum has three consequences. 
First, we can define the index of the active chart as
\begin{equation}
    \label{eq:learned selector}
    \hat m(\Dl,\tau) = \underset{m}{\text{argmin}}\, N\bigl(\phi(\bm{z}_m(\Dl)),\ t\bigr) \,\,,
\end{equation} 
which is the network's guess of the nearest image. Clearly, we expect all charts to tie where $V^{*}$ is horizon-limited, while $\bm{s}_m = \bm{0}$ wherever $V^{*}$ is distance-limited (i.e. $\hat m=0$). This is shown explicitly in Figure \ref{fig:selector}, where we plot the chart selector $m(\Dl,\tau)$ for the geometry $\TT$ on the domain $[0,2)^2$. In this way, the pursuer is at the origin $(0,0)$ and its periodic-images lie at the other corners, and the cut loci \footnote{From now on, we include a tilde to distinguish the fundamental domain, $\tilde\Dl\in[-1,1)^n$, from the domain $\Dl\in [0,2)^2$ of Figure \ref{fig:selector}.} $\Tilde\Delta_i=\pm1$ appear as midlines $\Delta_i=1$.

Moreover, the minimum in \eqref{eq:min of charts NN} determines the continuity properties of $V$. In fact, if one crosses the faces $\Tilde\Delta_i = \pm L/2=\pm1$ (i.e., the lines $\Delta_i=1$ in Fig. \ref{fig:selector}), it can be shown that the set of the possible displacements of the pursuer's images is left unchanged for any network $N$. This implies that the value function is continuous across $\Tilde\Delta_i = \pm L/2=\pm1$ and the ridge appears where two or more charts tie. On the other hand, the hyperplanes $\Tilde\Delta_i = 0$ (i.e., $\Delta_i=0,2$ in Fig. \ref{fig:selector}) require separate consideration. This is because $\Tilde\Delta_i$ changes sign, and so the set of the possible displacements changes. The continuity across $\Tilde\Delta_i = 0$ is guaranteed only when the minimum is attained by a chart with $(s_m)_i = 0$, otherwise, it would require the network to take the same value at $(z_m)_i=-L$ and $(z_m)_i=+L$. Notice that the condition $\bm{s}_m=0$ is the same when the viscosity solution $V^*$ is distance-limited. 

Finally, a minimum of smooth functions is semiconcave, so it can form peaks but not valleys~\citep{cannarsa2004semiconcave}. The network therefore cannot form a valley between charts, although it can still approximate one inside a single chart.

\subsubsection{The geometric NN}
\label{sec:geometric}

The geometric NN receives the geodesic distance $d$ directly. Among the $2^{n}$
image displacements $\bm{z}_m$ of Eq.~\eqref{eq: s_m image displacements} it
first selects the shortest one,
\begin{equation}
    \label{eq:nearest image}
    m^{*}(\Dl) = \underset{m}{\operatorname{argmin}}\ \|\bm{z}_m\| \, ,
    \qquad
    \bm{z}_{*}(\Dl) := \bm{z}_{m^{*}(\Dl)} \, ,
    \qquad \|\bm{z}_{*}\| = d \,.
\end{equation}

Then the network is given as inputs the $\phi$--features of that displacement alone,
\begin{equation}
    \label{eq:geometric}
    V(\Dl,\tau) = N\bigl(\phi(\bm{z}_{*}(\Dl)),\ t\bigr),
    \qquad
    \phi(\bm{z}_{*}) = \bigl(\bm{z}_{*},\ \bm{z}_{*}^{2},\ \|\bm{z}_{*}\|\bigr) \,.
\end{equation}
In the fundamental domain, with the pursuer at the origin, the nearest image is the pursuer itself, so that $(z_*)_i = \Delta_i$ and the network simply receives $(\Dl, \Dl^{2}, d)$. Compared with Eq.~\eqref{eq:min of charts NN}, the minimum over $m$ has moved inside the network: the min-of-charts NN evaluates $N$ on all $2^{n}$ charts and keeps the smallest value, whereas the geometric NN keeps the shortest chart and evaluates $N$ on it, so it needs one network evaluation per point and never has to learn which image is the nearest. The geometric NN is therefore the most informed about the geometry and the geodesic distance of all four configurations. However, nothing guarantees its periodicity or the sign of its kinks. Across the faces of the fundamental domain the component $(\bm{z}_{*})_i$ jumps from $+L/2$ to $-L/2$, while $(\bm{z}_{*})_i^{2}$ and $\|\bm{z}_{*}\|$ are continuous. This implies that $V$ is continuous on the torus only if $N$ takes the same value at $(\bm{z}_{*})_i = \pm L/2$, which we measure after training. Finally, we note that nothing prevents $N$ from forming a fold.

\section{Results}
\label{sec:results}
In this section, we present and compare the results obtained for every run, analyzing the successes and failures of the different NN configurations for the different geometries. We present the following quantities to diagnose the nature of the solutions:
\begin{enumerate}
    \item The best validation loss reached during training, namely the minimum of \eqref{eq:loss} over the run on a fixed validation set (for the soft-periodic NN we also include the additional periodicity term \eqref{eq:seam}).
    \item The \emph{interior} mean squared residual (MSR)
    \begin{equation}
        \mathrm{MSR} = \dfrac{1}{P}\sum_{P} r[V_\text{NN}]^{2}\,\,,
    \end{equation}
    at $P$ points in $\{d > 3\varepsilon,\, \, 0.5 < \tau < 7.5\}$, away from the terminal and capture boundary layers. 
    \item The mean gradient norm $\|\nabla V\|$ in the capture band $d \in (\varepsilon, 1.5\varepsilon]$, to be compared with $1/\dv = 5$, which is its gradient on the viscosity solution wherever the band is distance-limited (everywhere for $\tau > 0.2$).
    \item The root mean square (RMS) against the analytic solution $V^{*}$: \begin{equation}
        \mathrm{RMS} = \sqrt{\frac{1}{N_s}\sum_{i=1}^{N_s} \Bigl[V_\text{NN}(\Dl_i, \tau_i) - V^*(\Dl_i, \tau_i)\Bigr]^{2}}\, \, ,
        \label{eq: rms}
    \end{equation}
    on $N_s=2\cdot 10^{5}$ states drawn uniformly with $\Delta_i \in [-0.98, 0.98]$, $d > \varepsilon$ and $\tau \sim U(0.05, 8)$.
    
    \item The share of the $N_s$ uniformly sampled states on which the value is correct to within $5\%$, in the sense $|V -V^{*}| <0.05\,\max(V^{*}, 1)$: a relative tolerance where $V^{*} \ge 1$ and an absolute one of $0.05$ below it \footnote{Since $V^{*}\to0$ as one approaches the capture radius $d(\Dl)\to\varepsilon$ and as $\tau\to0$, where it has a corner no smooth network reproduces exactly, a purely relative tolerance fails there for every network; the floor ($0.05$ in time units, $\varepsilon/8$ in position) removes that thin band. Without it the correct configurations read $96$--$98\%$ on $S^{1}$ and $99$--$100\%$ on $\mathbb{T}^{2}$, with no change in the conclusions.}.
    
    \item Additionally, in \ref{app:cert} we present a concavity test that can be used as a certificate to validate the reached solutions.
\end{enumerate}

The first three diagnostics and the concavity test do not require knowledge of the analytic solution $V^*$, whereas the fourth and fifth do.

\subsection{Overview}
\label{sec:results-overview}

\begin{figure}[t!]
\centering
\includegraphics[width=\textwidth]{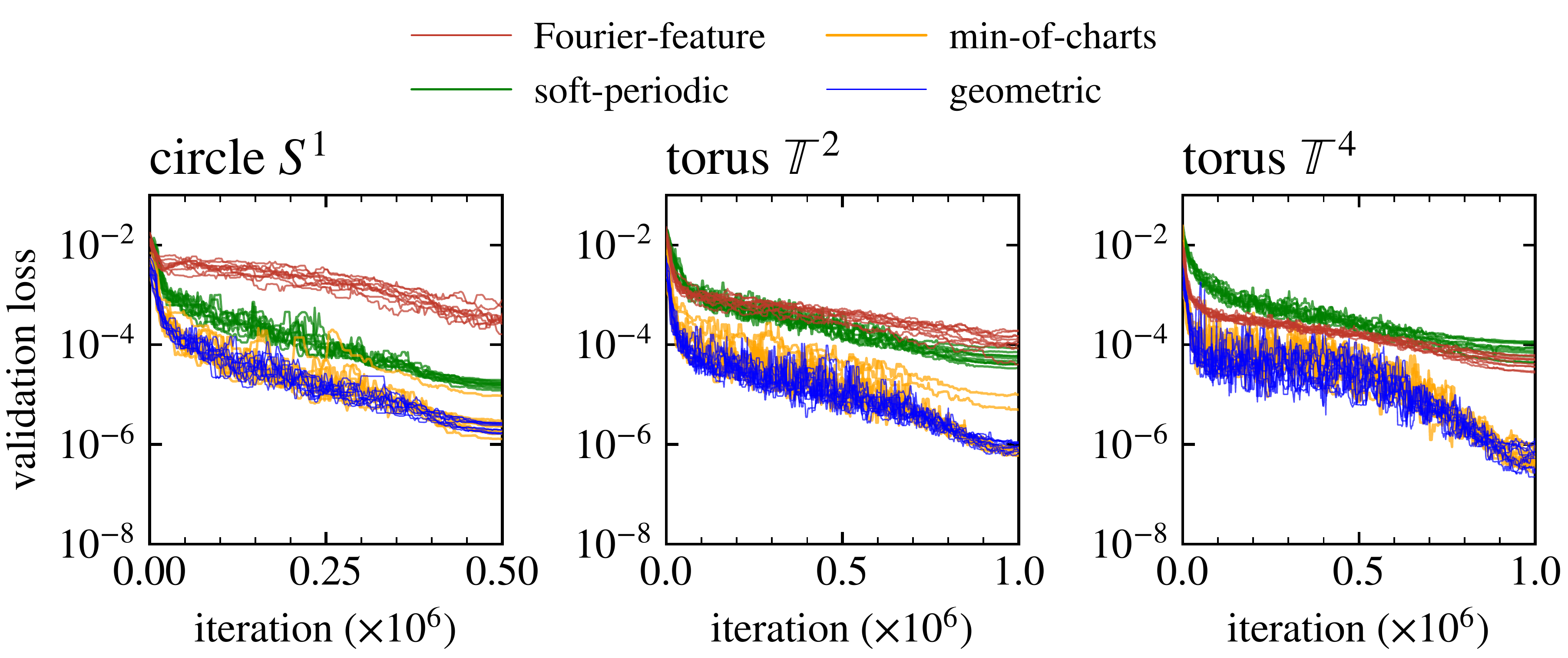}
\caption{Validation loss on $S^{1}$, $\TT$ and $\TF$ (left to right) for ten seeds per configuration, computed from Eq. \eqref{eq:loss} evaluated on a fixed random sampling.  The Fourier-feature NN (red) ends $70$--$180$ times above the min-of-charts (orange) and geometric (blue) NNs, which descend together. The soft-periodic NN (green), whose loss includes the additional term \eqref{eq:seam}, ends between the two groups on $S^{1}$ and $\TT$ and above the Fourier-feature NN on $\TF$.  All losses are still decreasing at the end of the budget.}
\label{fig:losses}
\end{figure}

Table~\ref{tab:campaign} lists the diagnostics of all $120$ runs, while Figure~\ref{fig:losses} presents the validation losses. At the residual level, all the results point in the direction of well-behaved, converged solutions. Every configuration reaches a small validation loss, and all but the soft-periodic NN on $\TF$ resolve the capture boundary layer to within $5\%$ of the analytic gradient (fourth column of Table~\ref{tab:campaign}). The min-of-charts and geometric NNs reach losses $70$--$180$ times below the Fourier-feature NN. However, comparing with the true viscosity solution $V^*$, we observe from Table \ref{tab:campaign} that the Fourier-feature NN is wrong over most of the domain on every geometry and for every seed. This implies that it converged to a folded solution of \eqref{eq:hji}, correct almost everywhere in residual but far from the viscosity solution in value. These two cannot be distinguished at the level of the loss function or the value of $\|\nabla V\|$ at capture (fourth column of Table~\ref{tab:campaign}). Along the same line, we see that the soft-periodic NN gives good results on $S^{1}$ and $\TT$, but then similarly fails on $\TF$. On the other hand, we can see that the two configurations whose input contains the geodesic distance, namely the min-of-charts and geometric NNs, find the correct result on most of the runs, apart from the three exceptions in the min-of-charts case, which we discuss further in Section~\ref{subsec: falure mode analysis}.

\begin{table}[!tbp]
\centering
\scriptsize
\setlength{\tabcolsep}{3pt}
\begin{tabular}{lcccc}
\toprule
Configuration & Best validation loss & RMS error against $V^{*}$ & $\|\nabla V\|$ at capture & Within $5\%$ of $V^{*}$ \\
\midrule
\multicolumn{5}{l}{\emph{\textbf{Circle $S^{1}$}}} \\
Fourier-feature & $1.5\ [0.8,4.6]\times10^{-4}$ & $1.36\ [0.93,1.86]$ & $5.00\ [4.99,5.01]$ & $44\ [19,60]\%$ \\
soft-periodic & $1.5\ [1.2,1.9]\times10^{-5}$ & $8.7\ [7.4,9.5]\times10^{-3}$ & $5.00\ [5.00,5.01]$ & $100\%$ \\
min-of-charts & $2.2\ [1.2,9.0]\times10^{-6}$ & $4.4\ [3.6,6.8]\times10^{-3}$ & $5.00\ [5.00,5.00]$ & $100\%$ \\
geometric & $\bm{2.1\ [1.6,2.6]\times10^{-6}}$ & $\bm{4.3\ [3.9,4.8]\times10^{-3}}$ & $5.00\ [5.00,5.00]$ & $100\%$ \\
\midrule
\multicolumn{5}{l}{\emph{\textbf{Torus $\mathbb{T}^{2}$}}} \\
Fourier-feature & $9.8\ [3.7,14.0]\times10^{-5}$ & $1.38\ [1.14,1.75]$ & $5.00\ [4.99,5.03]$ & $33\ [24,46]\%$ \\
soft-periodic & $4.8\ [3.2,8.2]\times10^{-5}$ & $2.0\ [1.9,2.5]\times10^{-2}$ & $5.00\ [4.99,5.01]$ & $100\ [99,100]\%$ \\
min-of-charts & $\bm{6.5\ [5.0,79.2]\times10^{-7}}$ & $\bm{3.0\ [2.6,12.0]\times10^{-3}}$ & $5.00\ [4.99,5.02]$ & $100\%$ \\
geometric & $7.5\ [5.2,9.7]\times10^{-7}$ & $3.3\ [2.9,4.1]\times10^{-3}$ & $5.00\ [4.99,5.00]$ & $100\%$ \\
\midrule
\multicolumn{5}{l}{\emph{\textbf{Torus $\mathbb{T}^{4}$}}} \\
Fourier-feature & $4.3\ [2.7,5.9]\times10^{-5}$ & $2.61\ [2.33,3.16]$ & $5.15\ [4.78,5.40]$ & $16\ [9,19]\%$ \\
soft-periodic & $7.1\ [4.3,11.3]\times10^{-5}$ & $5.1\ [1.9,11.2]\times10^{-1}$ & $2.99\ [2.75,3.71]$ & $69\ [47,83]\%$ \\
min-of-charts & $\bm{2.4\ [1.3,6.4]\times10^{-7}}$ & $\bm{8.7\ [6.4,98.8]\times10^{-3}}$ & $5.00\ [4.92,5.17]$ & $100\ [93,100]\%$ \\
geometric & $2.9\ [1.3,4.8]\times10^{-7}$ & $8.9\ [6.6,15.6]\times10^{-3}$ & $5.03\ [4.94,5.12]$ & $100\%$ \\
\bottomrule
\end{tabular}
\caption{Metrics for the four configurations on the three geometries $\mathbb{T}^{n}$ ($n = 1, 2, 4$), ten seeds each, with the same training budget within a geometry. ($5\times10^{5}$ iterations on $S^{1}$, $10^{6}$ on $\TT$ and $\TF$).  Entries are the median over seeds with the seed range in brackets (in the units of the entry).  Best validation loss is Eq. \eqref{eq:loss} on a fixed validation set; for the soft-periodic configuration it includes the extra term of Eq.~\eqref{eq:seam}.  The RMS is computed against the analytic $V^{*}$ on $2\times10^{5}$ uniform states ($\tau \sim U(0.05, 8)$, $d > \varepsilon$); $\|\nabla V\|$ at capture is the mean gradient norm in the band close to the capture set, $d \in (\varepsilon, 1.5\varepsilon]$, $\tau \sim U(0.5, 7.5)$, whose analytic value is $5$; the last column is the share of the same $2\times10^{5}$ states within $5\%$ relative error of $V^{*}$ (with the caveat in the item 5 at the beginning of the section). Bold marks the best median within each geometry (lowest loss and error, highest share within $5\%$), with ties all marked.}
\label{tab:campaign}
\end{table}

\begin{figure}[htp!]
\centering
\includegraphics[width=\textwidth]{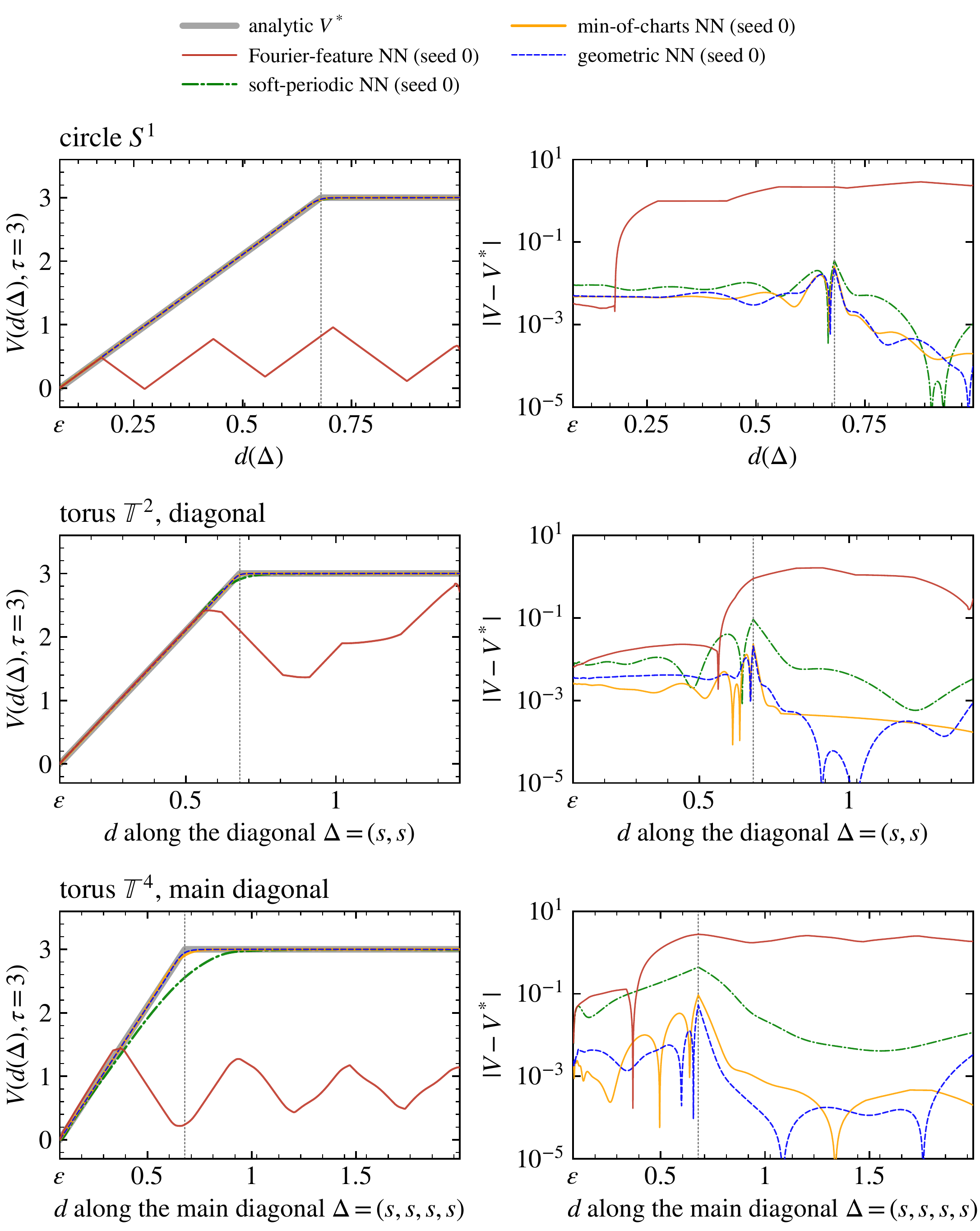}
\caption{Value function (left) and absolute error (right) of every configuration, against the analytic $V^{*}$ (gray),  evaluated at $\tau = 3$ of an example seed (seed~0). We show $S^{1}$ (top), a cut along the diagonal of $\TT$ (middle) and a cut along the main diagonal of $\TF$ (bottom). The dotted line is the cone $d_c = \varepsilon + \dv\,\tau$, beyond which $V^{*} = \tau$.}
\label{fig:profiles}
\end{figure}

In Figure \ref{fig:profiles}, we present the learned solutions $V(\Dl,\tau)$ for the different NN configurations, taking a single time slice $\tau=3$. In particular, a single spatial slice is shown along the diagonal for $\TT$ and along the principal diagonal for $\TF$.
From the Fourier-feature NN on $S^1$ of Figure \ref{fig:profiles}, notice that the solution shows a characteristic zig-zag feature. This is because the network has converged to a folded, almost-everywhere spurious solution, which satisfies the residual in Eq. \eqref{eq:hji} everywhere except in a measure-zero set of points, but differs wildly in value from the viscosity solution $V^*$. In fact, one can see how $\|\nabla V\|\simeq 1/\Delta v$, which is the quantity that the residual itself is constraining. Spurious solutions appear also in the $\TT$ and $\TF$ geometries. However, since in these cases one has to choose a specific direction in which to cut the tori, the zig-zag feature gets slightly deformed. Hence, the Fourier-feature NN converges to a wrong folded solution on every seed and every geometry studied; the mechanism is further discussed in Sec. \ref{subsec: falure mode analysis}.

\begin{figure}[htp!]
\centering
\includegraphics[width=\textwidth]{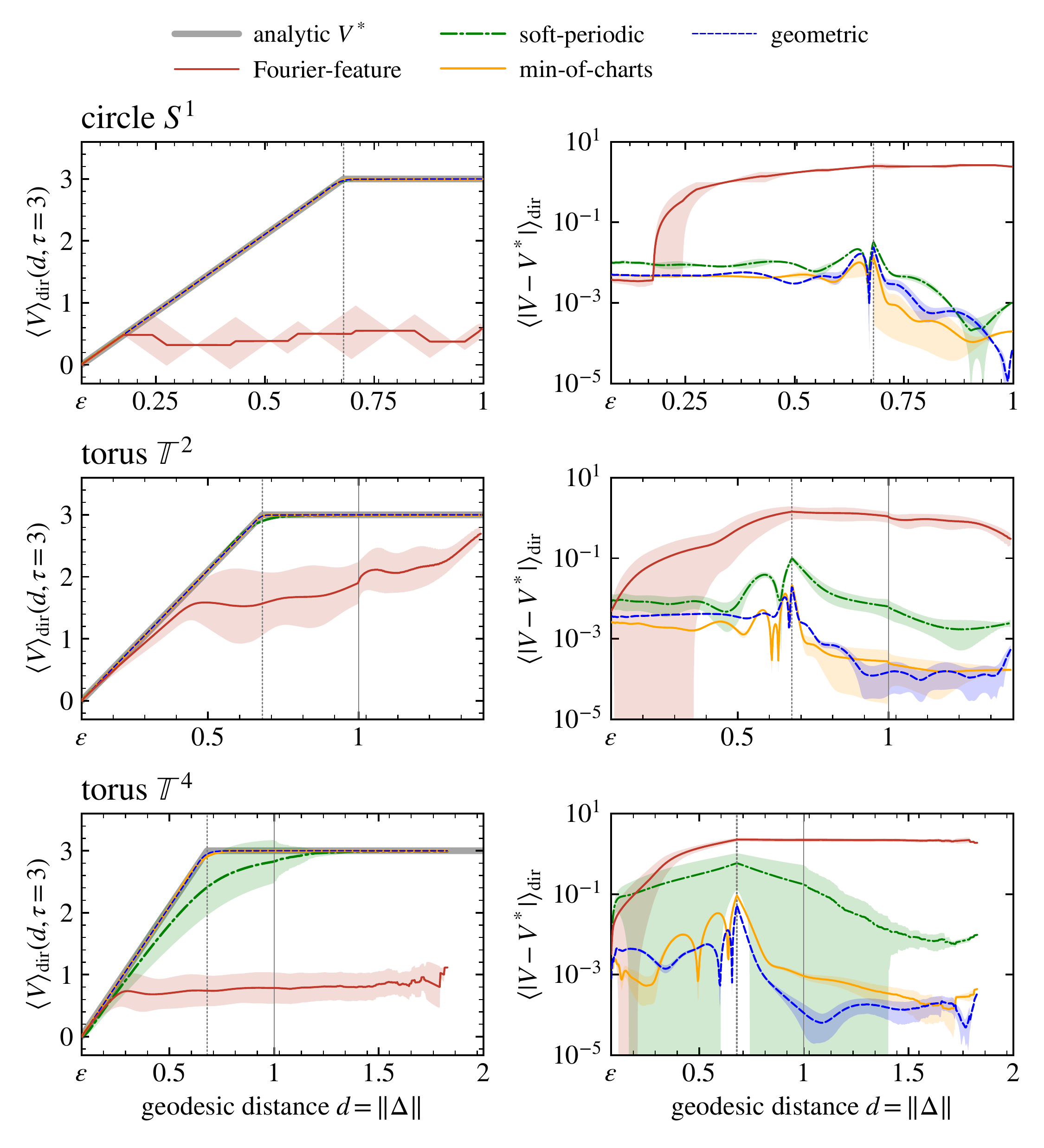}
\caption{Value function $V$ averaged over directions: at each geodesic distance $d$, $V$ is evaluated on all directions of the sphere $\|\Dl\| = d$ inside the fundamental domain ($2$ on $S^{1}$, $720$ on $\TT$, $512$ on $\TF$), and the mean and standard deviation are computed. The solid vertical line marks $d = 1$, beyond which only part of the sphere lies inside the domain. The dotted vertical line marks the cone, $d_c$. Same run and panel disposition as Fig. \ref{fig:profiles}. \textbf{Left:} mean over directions with a $\pm1\sigma$ band. The analytic $V^{*}$ depends on $d$ alone (not on direction), so any spread in the $1\sigma$ band is error due to angular variation of the learned solutions. \textbf{Right}: mean absolute error with its $\pm1\sigma$ band.}
\label{fig:angavg}
\end{figure}

On the other hand, we see that the remaining three parameterizations (the soft-periodic, the min-of-charts and the geometric NNs) converge to good approximations to the viscosity solution (see left plots in Fig. \ref{fig:profiles}). Concretely, for the soft-periodic case, the solution obtained for $\TT$ and $\TF$ is visually less converged; the absolute error in the right plots of Fig. \ref{fig:profiles} shows how this parameterization decreases in accuracy as we increase the dimension of the torus $\mathbb{T}^n$ from $n=1$ to $n=4$. We check that this feature is dependent on the direction in which we cut the tori in $n=2$ and $n=4$ (see Fig. \ref{fig:angavg}). 

We report in Figure \ref{fig:angavg} the mean and standard deviation values of the NN solutions averaged over different directions for a fixed distance $d=\|\Dl\|$, for the different values of $d$. A widening of the $1\sigma$ band reflects higher variability of the solutions at that specific distance when changing the direction of the spatial cut. From Figure \ref{fig:angavg}, we see that the soft-periodic band increases with the dimension of the torus, indicating that its accuracy depends increasingly on the spatial direction\footnote{The diagonal directions depicted in Fig. \ref{fig:profiles} are the ones with the smallest error.}. We discuss this feature in the analysis of the soft-periodic failure modes in Section \ref{subsec: falure mode analysis}. The tiny bands in Figure \ref{fig:angavg} of both the min-of-charts NN and the geometric NNs reflect an almost-null variability with direction. 

\begin{figure}[H]
\centering
\includegraphics[width=\textwidth]{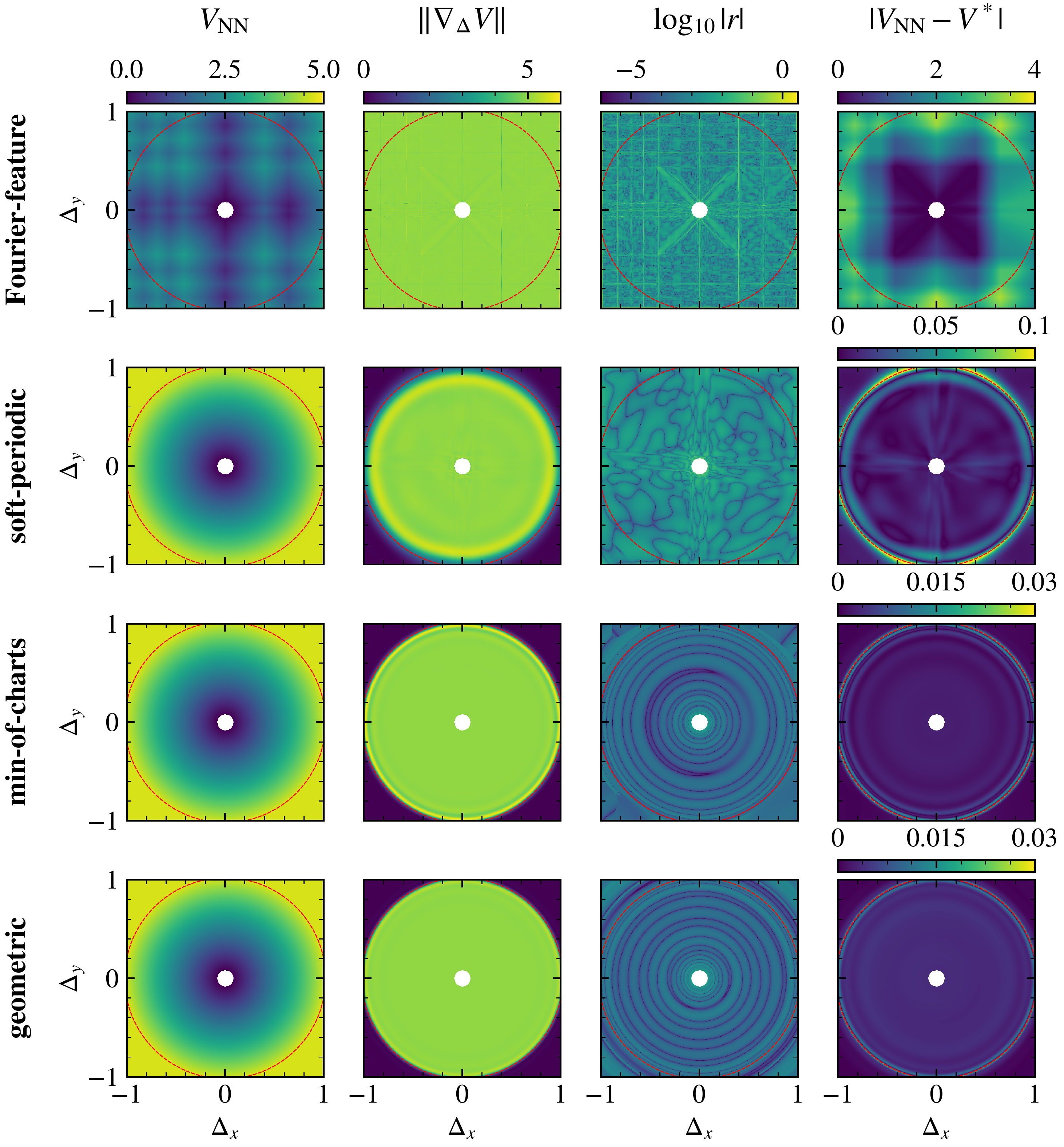}
\caption{Structure maps on $\TT$ at $\tau = 4.7$ (seed~0): from left to right columns, learned value $V_\text{NN}$, gradient norm $\|\nabla_{\Dl} V\|$, pointwise residual $|r|$ from Eq. \eqref{eq: residual HJI} (log scale) and absolute error $|V_\text{NN}-V^*|$, as a function of the wrapped coordinate $\Dl$. In each row, we show the four configurations. The dashed red circle is the cone $d_c$. The Fourier-feature NN has a gradient $\|\nabla_{\Dl}V \|$ of correct magnitude $1/\dv$ almost everywhere, but a lattice of folds instead of a cone in the value $V$, with the residual concentrated on the fold lines and an error that grows from the capture set outward reaching $O(1)$ values. The soft-periodic NN has the correct shape of the value $V$, with maximum absolute error of $10^{-1}$ concentrated around the cone. The min-of-charts and geometric NNs are within $3\cdot 10^{-2}$ of $V^*$ everywhere.}
\label{fig:maps2}
\end{figure}

\begin{figure}[H]
\centering
\includegraphics[width=\textwidth]{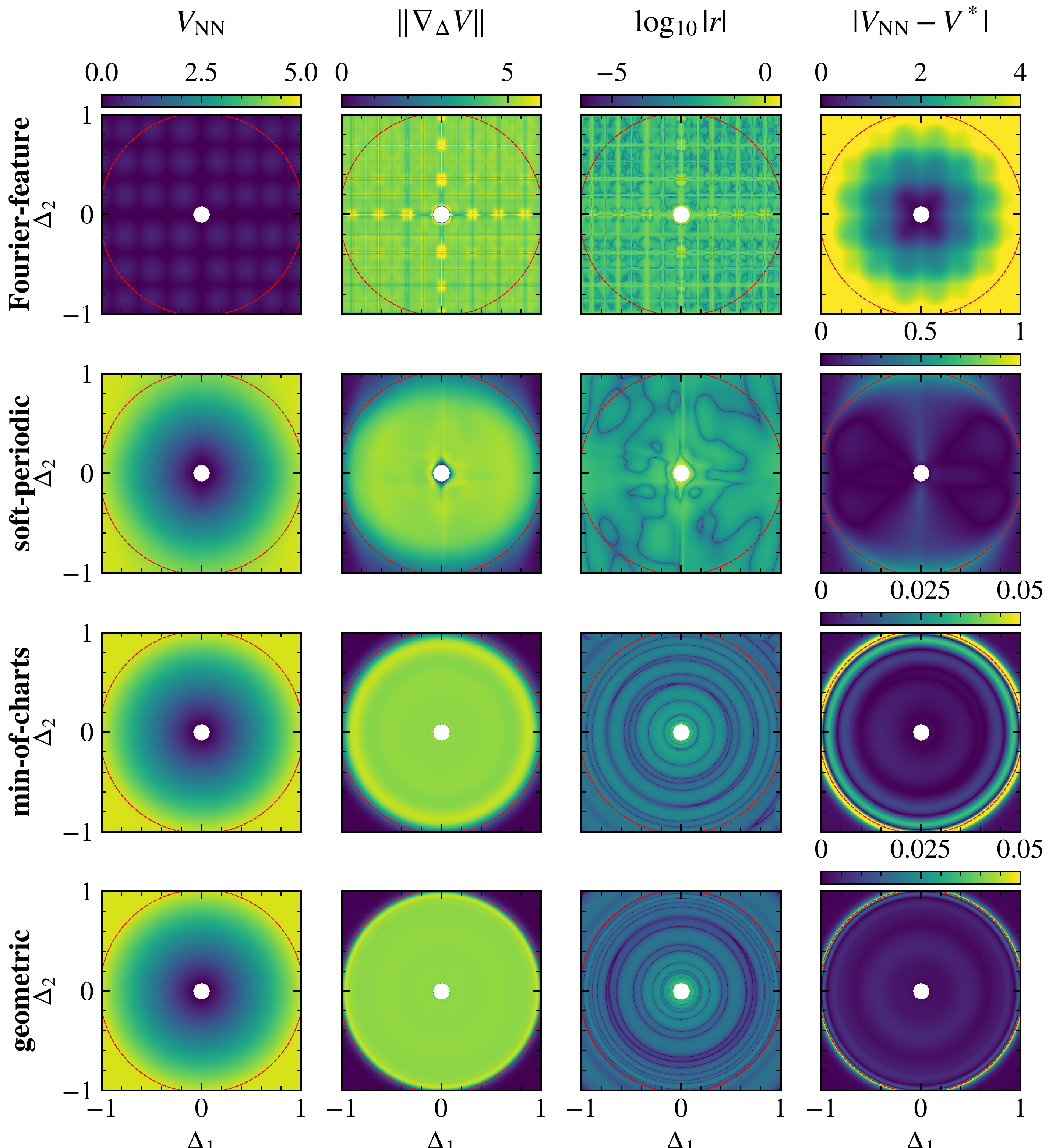}
\caption{Equivalent structure maps as in Fig. \ref{fig:maps2}, but for $\TF$, on the slice $\Delta_3 = \Delta_4 = 0$ at $\tau = 4.7$ (seed~0). The Fourier-feature NN shows the same fold lattice as Fig. \ref{fig:maps2}, with again $O(1)$ absolute errors over most of the slice, increasing as one moves away from the capture set. The soft-periodic NN has the right overall shape, but errors between 0.5 and 1, largest around the cone, especially along the line $\Delta_1 = 0$ and, more weakly, along part of the line $\Delta_2 = 0$. The min-of-charts and geometric NNs are within about $5\cdot 10^{-2}$ of $V^{*}$ on the whole slice, with the largest errors along the cone, and much lower values inside. In particular, notice that the first two configurations lack to recover the angular symmetry.}
\label{fig:maps4}
\end{figure}

These same points can be seen more explicitly in Figures \ref{fig:maps2} and \ref{fig:maps4}, in which we show a single $\tau$-slice of $V(\Dl,\tau)$ as a function of $(\Delta_x,\Delta_y)$ in the case of $\TT$ (Fig. \ref{fig:maps2}), and as a function of the first two coordinates $(\Delta_1,\Delta_2)$ on $\TF$, while keeping $\Delta_3,\Delta_4$ fixed (Fig. \ref{fig:maps4}). Notice that, in the third and fourth columns, the min-of-charts and geometric solutions have an almost-exact angular symmetry. This angular symmetry (exact in $V^*$) comes from the fact that the viscosity solution depends only on the (shortest) separation between the pursuer and evader and not on direction, i.e. the geodesic distance $d$. The min-of-charts and geometric NNs are informed on the geodesic distance through the possible nearest-image displacements $\bm{z}_m=\Dl-\bm{s}_m$, given as input features $\phi(\bm{z}_m)=(\bm{z}_m, \bm{z}_m^2, \|\bm{z}_m\|)$ (see Sec. \ref{sec:charts} and \ref{sec:geometric})\footnote{In these sections, we explain how the geodesic distance is given explicitly to the geometric NN via $\|\bm{z}_*\|=d$ of Eq. \eqref{eq:nearest image}, whereas the min-of-charts NN is given possible nearest-image distances $\|\bm{z}_m\|$ with $m=0,\dots, 2^n-1$ (see Eq. \ref{eq: s_m image displacements}), and has to select which of them is the true geodesic distance upon which to build the solution.}. Therefore, giving these geometry-informed inputs facilitates learning a solution with the relevant angular symmetry (specifically, through the feature $\|\bm{z}_m\|$). On the other hand, the soft-periodic configuration must learn this symmetry from the raw input coordinate $\Dl$, and the NN has to internally construct the modulus of the distance and learn that the solution only depends on it. The lack of accuracy in this symmetry for the soft-periodic case is most explicit in the third row of Fig. \ref{fig:maps4}.

Looking at the first row of Figs. \ref{fig:maps2} and \ref{fig:maps4}, we notice that the folded solution is far in value from the correct $V^*$, even if the value $\|\nabla_{\Dl} V\|$ matches the analytic value inside the cone, and the pointwise residual $|r|$ is small everywhere, but on a set of points of measure-zero \footnote{The set of points of measure-zero correspond to the convex kinks in the folded solutions, seen explicitly in the periodic, lattice-like residual of Fig. \ref{fig:maps2}, third column of the Fourier-feature NN.}.
This stresses again the fact that the HJI residual \eqref{eq: residual HJI} alone cannot constrain the NN-configurations enough to converge to the viscosity solution.

\subsection{Successes and failure-modes for the different NN configurations}
\label{subsec: falure mode analysis}

As already discussed in the previous section, the Fourier-feature NN fails by consistently converging to a folded, almost-everywhere spurious solution. The soft-periodic NN recovers $V^{*}$ on $S^{1}$ and $\TT$ in all $20$ runs and never gets a folded solution, but it fails on $\TF$. It is the least accurate of the three geometry-informed configurations, having a dependence of direction that increases with the dimension of the torus. The min-of-charts NN recovers $V^*$ in 27 of the total 30 runs, with median value RMS $0.0044$, $0.0030$, and $0.0087$ on $S^1$, $\TT$, and $\TF$, respectively. All sampled states are within $5\%$ of $V^*$ in the accurate runs, while the harmful $\TF$ run reaches $93\%$. The reason for the 3 failures is discussed below. Finally, the geometric NN recovers $V^{*}$ in all $30$ runs: value RMS $[3.9, 4.8]\cdot 10^{-3}$ on $S^{1}$, $[2.9, 4.1]\cdot 10^{-3}$ on $\TT$ and $[6.6, 15.6]\cdot 10^{-3}$ on $\TF$, every state within $5\%$ on every run, a resolved capture layer ($\|\nabla V\| = 4.94$--$5.12$), and the largest certificate ratio of any of its runs at $8.4$, below every folded and relabelled run (see \ref{app:cert}).

Let us analyze one by one the successes and failure modes of each configuration, and the reasons for them. \\

\begin{figure}[t!]
    \centering
    \includegraphics[width=\linewidth]{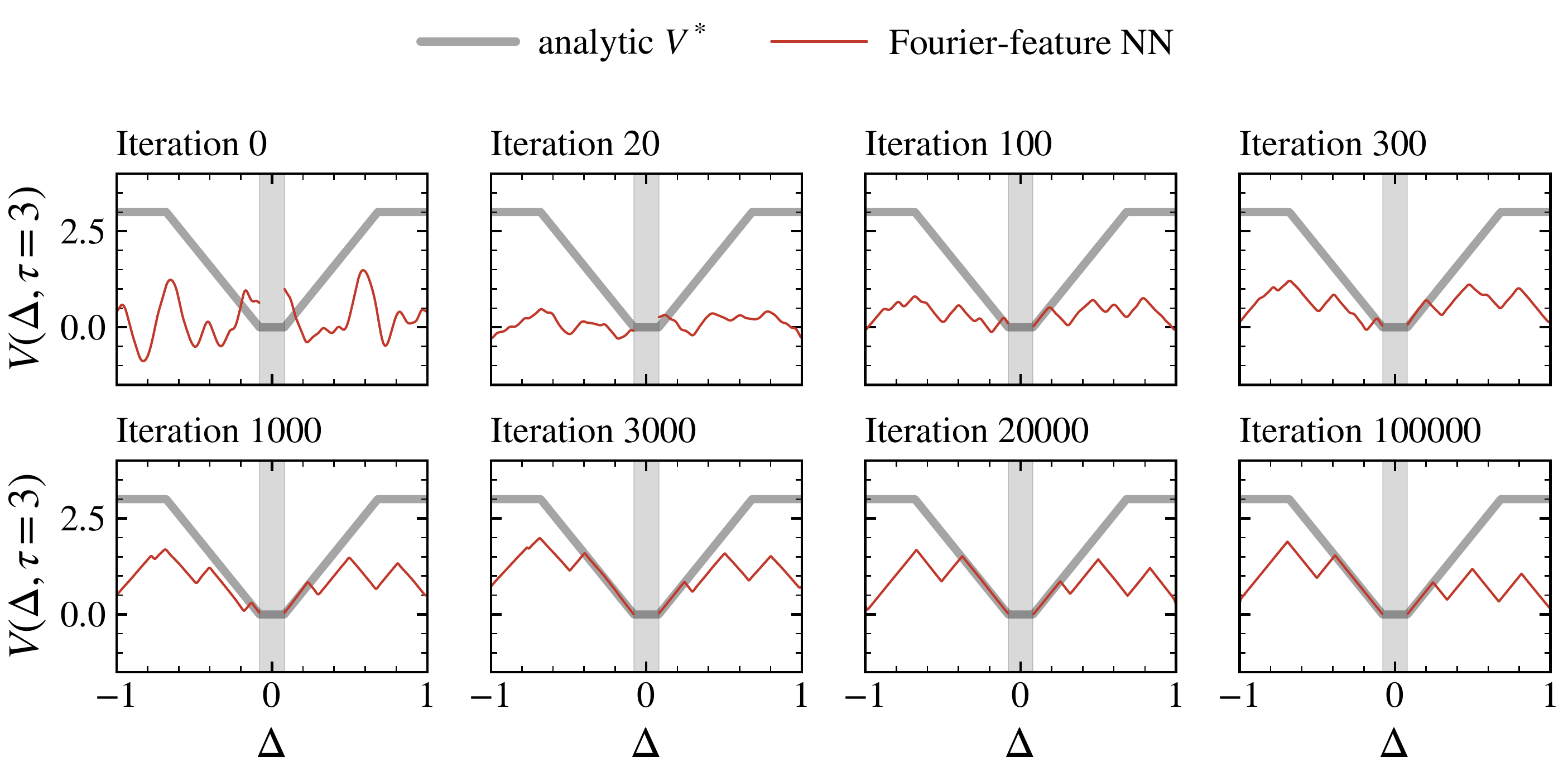}
    \caption{Evolution of Fourier-feature NN into almost-everywhere, folded solutions. The value function $V(\Delta, \tau = 3)$ of the Fourier-feature NN on $S^{1}$ (seed~0) is shown as a red line at different stages of a $10^{5}$--iteration run, from the initialization. The analytic $V^{*}$ is shown as a gray line. The gray band is the capture set $\|\Delta\| \le \varepsilon$.}
    \label{fig:fold_evolution}
\end{figure}

\noindent\textbf{Fourier-feature NN and folded solutions.} \\
Figure \ref{fig:fold_evolution} shows snapshots of a shorter run of $10^{5}$ iterations of the Fourier-feature NN on $S^{1}$ (seed~0), illustrating how a folded solution is reached during training. The mechanism is the one of Section~\ref{sec:analysis}: the HJI residual is a function of $|\partial_\Delta V|$ and $\partial_\tau V$, so it is minimised equally well by a slope of either sign. The capture condition fixes the value on the capture set, and the terminal condition fixes the level of the plateau at $V(\Dl,0)$, but not where the ramp should reach it. 
Starting from a small, oscillating, random function built from Fourier modes, the network can reduce the residual by giving the ramp the right magnitude wherever it is, with whatever sign it already has. Once a fold has formed, straightening it costs an $O(1)$ residual over a region. The training thus favours sharpening the folds rather than removing them. Although the loss continues to decrease, nothing in the curves suggests that longer training would make the network escape these folded solutions. Hence, despite the small residuals, all $30$ Fourier-feature runs end close to one of the non-viscosity almost-everywhere solutions that the theory of generalised solutions rules out~\citep{subbotin1995}.

Residual-driven adaptive sampling~\citep{daw2023r3} is a natural strategy to investigate, since the residual on the fold lines is well above its background (Figures~\ref{fig:maps2} and~\ref{fig:maps4}, third column).  We did not test whether adaptive sampling removes the folds or instead causes the network to sharpen them further. On a related one-dimensional problem, the authors in \citet{kim2026monotone} found that refinement based on the continuous residual placed almost no additional points at the spurious kink and left the error unchanged, whereas refinement based on a finite-difference residual removed most of it. The fact that residual training can produce such solutions was known~\citep{esteve2024fdls, kim2026monotone, golpashin2026viscosity}. Our results in Table~\ref{tab:campaign} show that these folded solutions are the typical outcome. Although in all 30 runs performed with Fourier-features we obtain these failures, they are completely invisible at the loss level.\\

\noindent\textbf{Soft-periodic NN and angular and dimension dependence.} \\
The soft-periodic NN is told nothing but the periodic identifications. It is the least accurate of the four geometry-informed configurations, apart from the Fourier-feature NN, with an error $2$--$7$ times and a loss $7$--$70$ times that of the min-of-charts and geometric NNs, residual maps with structure that the other two lack (see Figures~\ref{fig:maps2} and~\ref{fig:maps4}), and a jump across the faces of $(2$--$7)\times10^{-4}$ on $S^{1}$ and $(1.5$--$2.8)\times10^{-3}$ on $\TT$. Moreover, it shows an error that depends on direction: on $\TT$ it is several times larger along the axes than along the diagonal (Figure~\ref{fig:angavg}, middle row), where the capture points, drawn as cube-normalized random vectors, are denser.

On $\TF$, the soft-periodic NN fails to recover $V^*$ in all runs. The RMS ranges $0.19$--$1.12$, and only $47$--$83\%$ of the states lie within $5\%$ of $V^{*}$. 
The gradient norm near the capture boundary is only $\|\nabla V\| = 2.75$--$3.71$, below the expected value $1/\dv = 5$, indicating that the capture boundary layer is not resolved. 
The mean value error on the distance-limited branch ranges from $-0.2$ to $-1.2$, hence the solution value is too low everywhere, with the deficit largest along the coordinate axes (Figure~\ref{fig:maps4}, second row). 
The reason for this is sampling. The network has to learn the radial structure of $V^{*}$, the crease at the capture set $d = \varepsilon$ and the gradient $1/\dv$ just outside of it. It has to learn these structures from the residuals alone, that see the near-capture region only at sampled interior points close to it. For a fixed batch of $2048$ uniform points, the fraction in the band $d \in (\varepsilon, 1.5\varepsilon]$ scales with $\varepsilon^{n}$: about $84$ points per batch on $S^{1}$, $13$ on $\TT$ and $0.1$ on $\TF$. Therefore, in four dimensions the band that fixes the gradient at the crease is rarely sampled. The $512$ capture points fix the value $V = 0$ on the sphere $d = \varepsilon$, but not the slope leaving it. Similarly the $512$ face pairs of the periodic condition \eqref{eq:seam} cover the four faces of $\TT$ with $64$ pairs per unit length, but in the case of $\TF$ cover the eight three-dimensional faces of it with $8$ pairs per unit volume. The result is that the network smooths the crease, with a gradient of $3$ instead of $5$. Since the HJI equation propagates the value outward along the rays from the capture set, a crease that is too flat lowers the ramp everywhere, which is what the maps show.\\

\noindent\textbf{Min-of-charts NN failures.}\\
The min-of-charts NN recovers $V^{*}$ in $27$ of its $30$ runs, and fails very moderately in 3 runs. The reason is as follows. On the distance-limited states, our runs show a learned selector $\hat m(\Dl) $ (see Eq. \ref{eq:learned selector}) with binary outcomes: in every run it picks either the identity image ($26$ runs, $100\%$ of the tested states) or a single wrapped image ($4$ runs, $\geq 99.6\%$ of the tested states). The explanation on why this happens follows from Section~\ref{sec:charts}: at initialization and during training, the shared network $N$ is free to be smallest on a wrapped image everywhere. However, a wrapped image is not the geodesic distance, which corresponds to $\bm{s}_m=\bm{0}$ (no wrapping) in the fundamental domain $(-1,1]^n$. If this happens, the NN then has to approximate $V^{*}$ from the features of that image (which is the ``wrong shortest path''), whose wrapped coordinate jumps from $-2$ to $+2$ across the hyperplane $\Delta_i = 0$, where $V^{*}$ itself is smooth. $V$ is continuous there only if $N$ takes the same value at both ends of the jump, which it learns only approximately. Hence, the error concentrates on that hyperplane and, beyond the cone, in the lobes along the wrapped axis.

How much this failure-mode affects the quality of the actual solution depends on where the jump sits. On $S^{1}$ the jump is at $\Delta = 0$, a single point inside the capture set, where the capture condition fixes the value anyway, so the ``failed'' run is as accurate as its ``healthy'' siblings (value RMS $0.0068$ against $0.0036$--$0.0051$). However, for $n\geq 2$ the locus $\Delta_i=0$ is a codimension-one set that crosses the whole domain and meets the capture set on a negligible part of itself. For example, imagine that on $\TT$ the input of N for the ``wrong'' selected chart is wrapping around $x$, i.e. $\bm{z} = \Dl+ (-2 \text{sgn} \Delta_x, 0)$. The image-distance's first coordinate is $\Delta_x-2$ for $\Delta_x > 0$ and $\Delta_x + 2$ for $\Delta_x < 0$, so it jumps from -2 to +2 across the whole line $\Delta_x = 0$. That line has length 2 inside the fundamental domain, and only the segment $\|\Dl\|\leq \varepsilon$ lies inside the capture set, i.e. $8\%$ of its length. On the remaining $92\%$ of the line, the true value $V^*$ is smooth and nonzero, but there is no boundary condition pinning it and the network has to make $N(\phi(\bm{z}))$ continuous across a jump of 4 in one of its inputs. The even features of the chart do not jump (the feature $z_1^2$ is 4 on both sides and $\|\bm{z}\|$ is continuous), so only the raw coordinate is discontinuous. Therefore, N learns to be approximately even in $\bm{z}$, being the error a quantification of how good the approximation is. On $\TF$ the hyperplane is a three-dimensional subtorus ${\mathbb{T}^3}$ represented in the fundamental domain by a cube of side $2$. Its intersection with the capture set is a three-dimensional ball of radius $\varepsilon$ which occupies a fraction 
$$
\frac{\frac{4\pi}{3}\varepsilon^3}{2^3}=\frac{\pi}{6}\varepsilon^3\simeq 2.7\times10^{-4}
$$
of this subtorus for $\varepsilon =0.08$.
On $\TT$ the two failed runs have value RMS $12\cdot 10^{-3}$ and $8.6 \cdot 10^{-3}$, $3$--$4$ times their healthy siblings, still with $100\%$ of the states within $5\%$. On $\TF$ the failed run has value RMS $98.8 \cdot 10^{-3}$, $5$--$15$ times its healthy siblings, with $93\%$ of the states within $5\%$, and an error ten times larger near the hyperplane than away from it (RMS $0.27$ for $|\Delta_2| < 0.05$ against $0.036$ for $|\Delta_2| > 0.3$). 
An internal analysis of the active chart (check of selector $\hat m(\Dl)$) flags all four runs, and the semiconcavity test of \ref{app:cert} flags only the three harmful ones (not the $S^1$ failure). 

The $V_{\rm min-of-charts}=\underset{m}{\text{min}}\,V_m$ therefore does not improve accuracy over the geometric NN in healthy runs (their value, residual and loss are essentially indistinguishable), and it adds a failure mode. The min-of-charts NN is, though, more expensive, requiring $2^{n}$ network evaluations per point instead of one.
Its main advantage is the sign guarantee of Section~\ref{sec:charts} on the kinks formed between charts, and also the fact that it is less informed on the geometry than the geometric NN. In the three harmful runs, the selector check identifies incorrect image choices, but the concavity guarantee did not prevent them because the jump of the selected chart's input lies outside what the guarantee covers. Giving $N$ only the even features $(\bm{z}^{2}, \|\bm{z}\|)$ would remove this jump at no extra cost. However, we have not retrained with this change. \\

\noindent\textbf{Geometric NN is consistently successful.} \\
The geometric NN recovers $V^{*}$ in all $30$ runs, with the metrics shown in Table \ref{tab:campaign} and detailed at the beginning of this subsection. Its only error is a narrow peak of up to about $0.1$ absolute error at the cone (Figures~\ref{fig:profiles} and~\ref{fig:angavg}), where it smooths the concave kink that neither it nor the min-of-charts NN receives explicitly. Smoothing that kink costs no residual (Section~\ref{sec:analysis}), so the peak is the expected price of a smooth network and does not spread.

This configuration is given the most direct representation of the spatial structure relevant to this problem, which is why it works better. Its input is the nearest-image displacement $\bm{z}_{*}$ of Eq.~\eqref{eq:nearest image} together with $(\bm{z}_{*}^{2}, \|\bm{z}_{*}\|)$, so the geodesic distance $d = \|\bm{z}_{*}\|$, the only spatial variable on which $V^{*}$ depends, is an input. This makes it much easier to create the crease at $d = \varepsilon$, the ridge at the faces, and the radial dependence, all of which the soft-periodic NN has to learn from samples that become too sparse in four dimensions. Unlike the min-of-charts NN, the geometric NN has a single chart, the nearest one, and no chart selector, so it avoids the wrong-image failure mode of the min-of-charts construction. In this case, nothing in Eq.~\eqref{eq:geometric} enforces periodicity or the sign of the kinks, so both have to be learned. The jump of $V$ across the faces, where $(\bm{z}_{*})_i$ passes from $+L/2$ to $-L/2$, probed around a face at distances $10^{-3}$ and $10^{-2}$ on $2000$ random states (with identical results), is $(0.4$--$4)\times10^{-3}$ on $S^{1}$, $(1.5$--$7.6)\times10^{-4}$ on $\TT$ and $(0.3$--$3.8)\times10^{-3}$ on $\TF$. This jump is finite but small against values of order one, so we see that the network learns to be nearly periodic.  Its $30$ runs are clean by outcome, not by design, and the semiconcavity test of \ref{app:cert} makes this checkable without ground truth.

\section{Conclusions}
\label{sec:conclusions}

In this work, we analyzed four physics-informed neural network (PINN) configurations for the value function of a pursuit-evasion game on three translation-invariant geometries, namely $S^{1}$, $\TT$ and $\TF$, with ten seeds each. We compared every trained network with the closed-form viscosity solution of the  Hamilton--Jacobi--Isaacs (HJI) equation on these geometries. We draw three conclusions.

First, a sampled residual does not select the viscosity solution of a first-order HJ equation. The standard Fourier-feature network converged in all $30$ runs to a folded almost-everywhere solution with a small and still decreasing loss, a resolved capture boundary layer and a value that is wrong (too low) over most of the domain. For this equation, almost-everywhere solutions cannot lie above the value; the folds cost the residual almost nothing, so the loss function is blind to them. We observe that the optimization forms these folds, since it is the cheapest way of transforming the initial state of the Fourier-feature NN into a (wrong) solution satisfying the HJI equation almost everywhere. Once the optimizer has formed the folds, it does not remove them.

Second, as the network receives more information about the geometry, it selects the viscosity solution more reliably. Out of $30$ runs for each configuration, the Fourier-feature, soft-periodic, min-of-charts and geometric networks recovered the value in $0\%$, $66.7\%$, $90\%$ and $100\%$ of them, respectively. The Fourier-feature network knows only that the domain is periodic, which it receives through smooth harmonics of the state. The soft-periodic network receives the wrapped state, so it can place the ridge on the faces, and periodicity is imposed through the loss. This configuration failed only on $\TF$. The min-of-charts network receives geometric information from the norms of the $2^{n}$ possible shortest distances to the evader from periodic-images of the pursuer (one of which is the true geodesic distance). The set of possible displacements are labeled as $\bm{z}_m$ (with $m=0,\dots, 2^n-1$), and takes the minimum over the network evaluated on each $\bm{z}_m$. The network has to select which candidate $\bm{z}_m$ is the actual shortest (geodesic) distance. Its three failures are runs in which it selected a wrong image, which the minimum does not prevent. The geometric network receives the geodesic distance directly, as the norm of the wrapped state $\Dl$, so it is also told which image is the nearest, and it needs one network evaluation per point instead of $2^{n}$. The geodesic distance carries both the ridge and the radial structure of the value function, and supplying it directly was enough for the network to recover the viscosity solution in every run.

Third, the dimension degrades what the network has to learn from samples, but not what it receives in its input. The band next to the capture set, where the slope of the value function $V(\Dl,\tau)$ is fixed, holds a fraction of the collocation points that scales as $\varepsilon^{n}$, about $84$ points per batch on $S^{1}$, $13$ on $\TT$ and $0.1$ on $\TF$. The soft-periodic network, which has to learn the radial structure from these points, recovered the slope at capture on $S^{1}$ and $\TT$ but not on $\TF$, where it reached a gradient of about $3$ instead of $1/\dv = 5$ and a value that is too low everywhere. The share of states on which the Fourier-feature network lies within $5\%$ of the value also falls with the dimension, from a median of $44\%$ on $S^{1}$ to $33\%$ on $\TT$ and $16\%$ on $\TF$. The min-of-charts and geometric networks, which receive the crease at capture and the radial structure through the geodesic distance, kept a median RMS with respect to the true solution $V^*$ below $10^{-2}$ on every geometry.

A wrong solution can also be detected without the value function, provided the Hamiltonian is convex and a bound on the curvature of the value is known. The semiconcavity test of \ref{app:cert} checks this bound from forward evaluations of the trained network. In our runs it flags the folded solutions and the three harmful min-of-charts runs, which the loss does not. 

Clearly, our study has some limitations. All value-level results are based on a single differential game on translation-invariant geometries, for which the value depends only on the displacement between the players and the  Hamiltonian is convex. The sign of the error of folded solutions and the semiconcavity test both rely on this convexity, so they do not apply to HJI equations whose Hamiltonian is neither convex nor concave. We fixed the network width at $128$, used $K=8$ Fourier harmonics, and kept the training budget fixed. We did not investigate whether alternative choices would change the behaviour of the Fourier-feature network. Moreover, all training losses were still decreasing at the end of the training budget. Our configurations also do not isolate the effects of the ridge structure and the radial information provided by the distance input. In addition, we exclude a neighbourhood of width $0.02$ around the faces of the fundamental domain from both residual and value-error evaluation, leaving the accuracy near the cut-locus ridge unassessed (\ref{app:training}).
Finally, we did not test whether alternative sampling strategies, such as residual-driven adaptive sampling~\citep{daw2023r3} or sampling points concentrated near the capture set, might remove the folds or improve the performance of the soft-periodic network on $\TF$.

The solution proposed here relies on prior knowledge of the relevant distance structure of the problem. For problems in which this structure or the singular set is unknown, loss-level selection~\citep{esteve2024fdls, bokanowski2026monotone} or time-marching curricula~\citep{bansaltomlin2021deepreach} are possible alternatives, but their ability to prevent or eliminate folded solutions remains an open question. When the Hamiltonian is convex and a curvature bound is known, a natural next step is to incorporate the semiconcavity test of \ref{app:cert} directly into a loss term that penalises second differences of the network above the curvature bound at randomly sampled points. 
This would push away training from folded solutions rather than relying on the diagnostics to detect them only afterwards.

 \section*{Acknowledgements}
{\emergencystretch=3em
P.T.-A.\ is supported by project PRE2022-102220, ``Dark Energy and the Origin
of the Universe'', funded by MCIN/\allowbreak AEI/\allowbreak 10.13039/\allowbreak 501100011033. Funding for the work of P.T.-A., P.T.-P., L.S., P.S.-V., N.T., M.T.\ and R.J.\ was partially provided by project PID2022-141125NB-I00 and the ``Center of Excellence María de Maeztu 2020-2023'' award to the ICCUB (CEX2019-000918-M), funded by MCIN/\allowbreak AEI/\allowbreak 10.13039/\allowbreak 501100011033. 
M.T. acknowledges support of the ``la Caixa'' Foundation fellowship (ID 100010434). The fellowship code is ``LCF/BQ/DI24/12070012''. The work of PS was supported by a grant from the ``la Caixa'' Foundation (ID 100010434), under the fellowship code LCF/BQ/DI23/11990074. \par}

\section*{Code availability}
\label{app:repro}
 
All code, loss traces and diagnostic outputs are available in the public GitHub repository \href{https://github.com/ptejerina/PENN}{PENN}. The trained networks are available upon request to the corresponding author.

\section*{AI use statement}

Claude Code was used to assist with coding, computations, and graphic generation. The ideas, inferences, and results all originate in the (human) brains of the authors. Claude and ChatGPT were used to check citations, spelling, grammar, and inconsistencies. The text was written by the authors, who are responsible for the results, findings and conclusions of this article.

\appendix
 
\section{Training details}
\label{app:training}

Interior points are drawn at every training step, with $\Dl$ uniform in \mbox{$[-L/2+\delta, L/2-\delta]^{n}$}, a margin $\delta = 0.02$ at the faces of the fundamental domain, and $\tau$ uniform in $(0,\tau_{\max})$.

They are kept if they lie outside the capture set, $d(\Dl) \ge \varepsilon$.

Terminal points are drawn on the slice $\tau = 0$ and capture points inside the capture set. The loss \eqref{eq:loss} uses $N_{\mathrm{int}} = 2048$ interior, $N_{0} = 512$ terminal and $N_{\varepsilon} = 512$ capture points per step, each term being a mean over its own batch, and $\lambda = 5$. The periodicity penalty \eqref{eq:seam} uses $N_{p} = 512$ pairs of points, redrawn at every step. Each pair is assigned a random axis, along which its two points have coordinates $\pm L/2$. The two points share the remaining coordinates and the time $t_l$, which is drawn as for the interior points. The penalty has the same weight $\lambda = 5$ as the boundary terms.

All configurations use the same network $N$, a multilayer perceptron with two hidden layers of width $128$ and $\tanh$ activations, and the same time features $(t,\ t^{2},\ \sin\pi t,\ \cos\pi t)$, with $t = \tau/\tau_{\max}$. The number of inputs is $16n + 4$ for the Fourier-feature NN, $n + 4$ for the soft-periodic NN, and $2n + 5$ for the min-of-charts and geometric NNs. One training point costs $2^{n}$ network evaluations for the min-of-charts NN and a single evaluation for the others. Therefore, even if batch sizes and iteration counts are the same for all configurations on a geometry, the cost per iteration is not.

Every configuration is trained on every geometry with
10 seeds, $120$ runs in total. Derivatives are computed exactly by automatic differentiation.  We use Adam~\citep{kingma2015adam} with a single cosine learning-rate schedule from $3\times10^{-3}$ to $5\times10^{-5}$, for $5\times10^{5}$ iterations on
$S^{1}$ and $10^{6}$ on $\TT$ and $\TF$.  The validation losses are computed via evaluation of the NNs on points drawn independently of the training batches, and displayed in Figure~\ref{fig:losses}.

The runs have been trained on Apple-silicon GPUs with PyTorch 2.14 and on NVIDIA L40S GPUs; training seeds are recorded for every run. A run took between one and $16$ hours, the $16$-chart min-of-charts runs on $\TF$ being the most expensive.

\section{The semiconcavity test as a solution certificate}
\label{app:cert}

From forward evaluations of the trained network alone, the semiconcavity test checks the one-sided curvature bound that the value function $V^{*}$ satisfies. Away from the capture set and the ridge, the spatial Hessian of $(d - \varepsilon)/\dv$ is $(\mathbb{I} - \hat{\bm{x}}\hat{\bm{x}}^{\top})/(\dv\,d)$, with $\hat{\bm{x}}$ the unit displacement from the nearest image. Its eigenvalues are $0$ and $1/(\dv\,d)$, and the minima with $\tau$ and over images in Eq. \eqref{eq: Vstar as min of min} add only concave kinks. Hence,
\begin{equation}
\lambda_{\max}\bigl(\nabla^{2} V^{*}(\Dl,\tau)\bigr)\ \le\
C^{*}(\Dl) := \frac{1}{\dv\, d(\Dl)}
\qquad\text{for } d(\Dl) > \varepsilon ,
\label{eq:cstar}
\end{equation}
and, in the second-difference form that also holds at the kinks, for every unit vector $\bm{u}$ and every $0 < h < d(\Dl) - \varepsilon$,
\begin{equation}
\frac{V^{*}(\Dl + h\bm{u},\tau) - 2V^{*}(\Dl,\tau) + V^{*}(\Dl - h\bm{u},\tau)}{h^{2}}
\ \le\ \frac{1}{\dv\,\bigl(d(\Dl) - h\bigr)} ,
\label{eq:secdiff}
\end{equation}
because $V^{*}$ lies below the smooth branch that attains the minimum at $\Dl$ and touches it there. A folded solution violates \eqref{eq:secdiff} next to every convex kink, where its second differences grow like $1/h$.

\paragraph{Procedure} At each of $4\times10^{4}$ probes drawn uniformly in the fundamental domain at $\tau = 4$ we compute the second difference $S_h(\Dl;\bm{u}_j)$ of the network along $M$ random unit directions $\bm{u}_j$ ($M = 8$; $M = 1$ on $S^{1}$) with $h = 0.02$, and $\hat\lambda_{\max}(\Dl) = \max_j S_h(\Dl;\bm{u}_j)$. Probes with $d \le 3\varepsilon$ are excluded, which keeps the right-hand side of \eqref{eq:secdiff} below $1.1\,C^{*}$, and so are probes next to a concave kink, detected as $\min_j S_h < -3\,C^{*}$, since admissible kinks are not evidence. A probe is flagged when $\hat\lambda_{\max} > 2\,C^{*}$, and for each run we report the flagged fraction and the largest ratio $\hat\lambda_{\max}/C^{*}$. The test costs $2M + 1$ network evaluations per probe and needs no derivatives, no mesh and nothing of $V^{*}$ beyond the bound $C^{*}$. The finite step is essential: an automatic-differentiation Hessian sees a kink formed by a minimum over charts only from its smooth side, whereas $S_h$ responds to it from any probe within $h$. A fold smoothed over a width $\ell$ still gives a ratio of order $d/\max(\ell, h)$, so smoothing hides it only in proportion to its width, while by Section~\ref{sec:analysis} the residual penalises a fold in proportion to its width: the two checks are complementary.

\begin{table}[!tbp]
\centering\footnotesize\setlength{\tabcolsep}{4pt}
\begin{tabular}{lccc}
\toprule
candidate ($\tau = 4$) & value RMS & flagged \% & max ratio $\hat\lambda_{\max}/C^{*}$ \\
\midrule
\multicolumn{4}{l}{\emph{circle $S^{1}$}} \\
$V^{*}$ (analytic) & --- & $0.0$ & $0.0$ \\
planted sawtooth & --- & $8.2$ & $100$ \\
Fourier-feature & $1.36\ [0.93,1.86]$ & $11.6\ [5.2,18.8]$ & $212\ [88,568]$ \\
soft-periodic & $8.7\ [7.4,9.5]\times10^{-3}$ & $\bm{0.0\ [0.0,0.0]}$ & $\bm{0.8\ [0.4,1.1]}$ \\
min-of-charts & $4.4\ [3.6,6.8]\times10^{-3}$ & $\bm{0.0\ [0.0,3.6]}$ & $1.8\ [0.4,5.9]$ \\
geometric & $\bm{4.3\ [3.9,4.8]\times10^{-3}}$ & $\bm{0.0\ [0.0,4.1]}$ & $1.5\ [0.6,8.4]$ \\
\midrule
\multicolumn{4}{l}{\emph{torus $\mathbb{T}^{2}$}} \\
$V^{*}$ (analytic) & --- & $0.0$ & $1.0$ \\
planted sawtooth & --- & $8.3$ & $135$ \\
Fourier-feature & $1.38\ [1.14,1.75]$ & $20.5\ [12.7,37.5]$ & $158\ [100,883]$ \\
soft-periodic & $2.0\ [1.9,2.5]\times10^{-2}$ & $0.8\ [0.2,2.7]$ & $4.3\ [2.4,7.4]$ \\
min-of-charts & $\bm{3.0\ [2.6,12.0]\times10^{-3}}$ & $1.3\ [0.0,4.1]$ & $2.6\ [1.2,24.2]$ \\
geometric & $3.3\ [2.9,4.1]\times10^{-3}$ & $\bm{0.0\ [0.0,3.4]}$ & $\bm{2.0\ [1.3,3.0]}$ \\
\midrule
\multicolumn{4}{l}{\emph{torus $\mathbb{T}^{4}$}} \\
$V^{*}$ (analytic) & --- & $0.0$ & $1.0$ \\
planted sawtooth & --- & $7.9$ & $167$ \\
Fourier-feature & $2.61\ [2.33,3.16]$ & $87.4\ [81.7,94.4]$ & $112\ [96,156]$ \\
soft-periodic & $5.1\ [1.9,11.2]\times10^{-1}$ & $2.9\ [0.9,9.3]$ & $18.1\ [4.5,73.4]$ \\
min-of-charts & $\bm{8.7\ [6.4,98.8]\times10^{-3}}$ & $\bm{0.0\ [0.0,1.5]}$ & $\bm{1.1\ [1.0,33.1]}$ \\
geometric & $8.9\ [6.6,15.6]\times10^{-3}$ & $\bm{0.0\ [0.0,0.5]}$ & $\bm{1.1\ [1.1,2.3]}$ \\
\bottomrule
\end{tabular}
\caption{Monte Carlo semiconcavity test at $\tau = 4$ ($4\times10^{4}$
probes; eight random directions each and one on $S^{1}$; $h = 0.02$; the region
$d \le 3\varepsilon$ and concave kinks are masked; a probe is flagged when
$\hat\lambda_{\max} > 2\,C^{*}$), applied to $V^{*}$, to a planted sawtooth (a
$\tau$-independent triangle wave in the first coordinate $\Delta_{x}$, of
period $0.5$ and amplitude $1.25$, so with slope of magnitude $1/\dv$ on its
flanks and convex kinks at its valleys, which shows what a clearly folded
function scores) and to all $120$ trained networks; median and range over the
ten seeds.  The value RMS column repeats Table~\ref{tab:campaign}
($\tau \in (0.05, 8)$, not only $\tau = 4$).  On $S^{1}$ the curvature of
$V^{*}$ vanishes away from its kinks, so its ratio is $0$ rather than $1$.
Bold marks the lowest median among the trained networks within each
geometry, with ties all marked.}
\label{tab:cert}
\end{table}

\begin{figure}[!h]
    \centering
    \includegraphics[width=\linewidth]{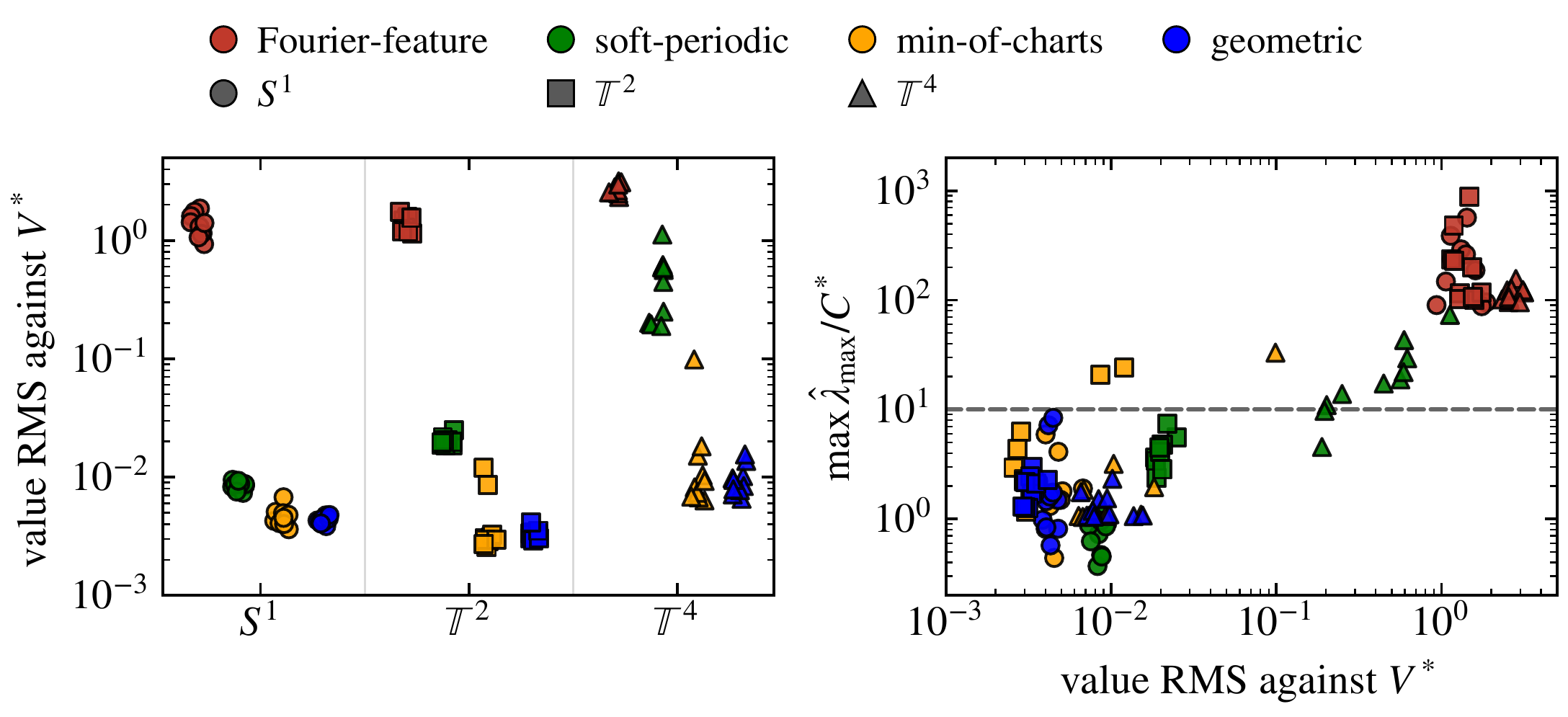}
    \caption{Seed census of the $120$ runs; one marker per run, colour by
    configuration, marker shape by topology.  \textbf{Left.} Value RMS against $V^{*}$ (Table~\ref{tab:campaign}) grouped by topology.  \textbf{Right.} The semiconcavity certificate $\max\,\hat\lambda_{\max}/C^{*}$ (Table~\ref{tab:cert}) against the same RMS; the dashed line is the ratio $10$.  Folded and relabelled runs lie above it and correct runs below it, while the $\mathbb{T}^{\,4}$ soft-periodic failures, smooth creases rather than folds, straddle it.}
    \label{fig: census}
\end{figure}

\paragraph{Findings} Table~\ref{tab:cert} and Figure~\ref{fig: census} show the following. $V^{*}$ has a maximum ratio of $1.0$ on the tori ($0.0$ on $S^{1}$, where its curvature vanishes) and no flagged probe; the planted sawtooth reaches $100$--$167$. All $30$ Fourier-feature runs reach $88$--$883$, with $5$--$94\%$ of the probes flagged, about two orders of magnitude above the correct runs, whose maximum ratios lie between $0.4$ and $8.4$ with at most $4.1\%$ of the probes flagged, at approximation ripples that the error maps show to be harmless. The three harmful min-of-charts runs reach $20.7$ and $24.2$ on $\TT$ and $33.1$ on $\TF$, against at most $6.3$ for their correct siblings, because of the near-jump of $V$ across the hyperplane where the input of their selected chart jumps (Section~\ref{subsec: falure mode analysis}); the harmless relabelled run on $S^{1}$, whose jump lies inside the capture set, scores $1.9$. We therefore judge a run by its maximum ratio and read a value above about $10$ as the sign of an inadmissible kink: the largest ratio among the correct runs is $8.4$ and the smallest among the folded and relabelled runs is $20.7$, whereas the flagged fraction alone would not separate them ($4.1\%$ for a correct run, $1.0$--$1.3\%$ for the relabelled runs on $\TT$). This threshold was chosen after seeing these results. The soft-periodic failures on $\TF$ are of a different kind, a crease that is too flat ($\|\nabla V\| \approx 3$) rather than a fold: eight of the ten runs still exceed $10$ (ratios $4.5$--$73$), but the two most accurate ones do not, and it is the capture gradient of Table~\ref{tab:campaign} that identifies them.

\paragraph{Scope} The test checks admissibility, not accuracy: it sees the sign of the kinks, not the level of the value, so it complements a residual or boundary check rather than replacing one; quantitative error bounds require monotone-scheme residuals~\citep{esteve2024fdls, bokanowski2026monotone} or a posteriori estimates~\citep{albert2002posteriori}. It needs a one-sided curvature bound for the solution class, available when the Hamiltonian is convex (or concave, with the sign reversed)~\citep{cannarsa2004semiconcave} but not for HJI equations whose Hamiltonian is neither~\citep[Remark~5.3.3]{cannarsa2004semiconcave}; here $C^{*}$ comes from the known structure of $V^{*}$, elsewhere it has to come from theory~\citep[Chapter~8]{cannarsa2004semiconcave}. Being Monte Carlo, the test can miss folds of very small extent, and a random direction $\bm{u}$ sees only the fraction $(\bm{u}\cdot\bm{\nu})^{2}$ of the curvature of a fold with normal $\bm{\nu}$ (with $M = 8$ on $\TF$ the best direction sees a median of $69\%$ of it), so the ratios reported for folds are conservative and $M$ has to grow with the dimension. The mask can also hide a fold within about $h$ of a concave kink; the folds of our trained networks are far coarser.

\end{document}